\documentclass[twoside]{article}
\usepackage{sibgrapi}
\usepackage{times}
\usepackage{amssymb}
\usepackage{amsmath}
\usepackage{fancyhdr}
\usepackage{subfigure}
\usepackage[dvipdfm]{hyperref}
\usepackage[dvipdfm]{graphicx}
\input{tcilatex}
\renewcommand{\headrulewidth}{0pt}
\begin{document}



\newcommand{\R}{\mathbb{R}}
\renewcommand{\v}{\mathbf{v}}
\renewcommand{\P}{\mathbf{P}}
\newcommand{\Q}{\mathbf{Q}}
\newcommand{\p}{\mathbf{p}}
\newcommand{\q}{\mathbf{q}}
\newcommand{\w}{\mathbf{w}}
\newcommand{\x}{\mathbf{x}}
\newcommand{\y}{\mathbf{y}}
\newcommand{\zb}{\overline{z}}
\renewcommand{\t}{\mathbf{t}}
\newcommand{\n}{\mathbf{n}}
\renewcommand{\b}{\mathbf{b}}
\newcommand{\T}{\mathbf{T}}
\renewcommand{\L}{\mathbf{L}}
\newcommand{\V}{\mathbf{V}}
\newcommand{\X}{\mathbf{X}}
\newcommand{\Y}{\mathbf{Y}}
\newcommand{\N}{\mathbf{N}}
\newcommand{\C}{\mathbf{C}}
\newcommand{\dl}{\Delta l}
\renewcommand{\k}{\kappa}
\newcommand{\norm}[1]{\Vert #1 \Vert}
\newcommand{\fsep}{\hspace*{\fill}}



\title{ Area distances of Convex Plane Curves and Improper Affine Spheres}
\shorttitle{Area Distances as Affine Spheres}

\author{
  Marcos Craizer\authortag{1},
  Moacyr Alvim\authortag{2}
  and
  Ralph Teixeira\authortag{3}
} 

\address{
  \authortag{1} Pontif\'{i}cia Universidade Cat\'{o}lica, \authortag{2} Funda\c{c}\~{a}o Getulio Vargas,  \\
   \authortag{3} Universidade Federal Fluminense -  Rio de Janeiro, Brazil. \\
   {\small \tt craizer@mat.puc-rio.br},\ \ {\small \tt moacyr.silva@fgv.br},\ \ {\small \tt ralph@mat.uff.br} \\
}

\abstract { The area distance to a convex plane curve is an
important concept in computer vision. In this paper we describe a
strong link between area distances and improper affine spheres.
This link makes possible a better understanding of both theories.
The concepts of the theory of affine spheres lead to a new
definition of an area distance on the outer part of a convex plane
arc. Also, based on the theory of discrete affine spheres, we
propose fast algorithms to compute the area distances. On the
other hand, area distances provide a good geometrical
understanding of improper affine spheres.
\\[2mm] \noindent
}

\maketitle

\section{Introduction}

The area distance of convex plane curves is an important concept
in computer vision. These distances can be useful in matching two
images of the same object obtained from different points of view
\cite{Giblin04}. It can also be seen as an erosion, a basic
concept of mathematical morphology (\cite{Moisan98}).

Improper affine spheres are surfaces in $R^3$ whose affine normals
at all points are parallel. In this paper, we point out the strong
connection between these area distances and improper affine
spheres. This connection is used in the development of the theory
of area distances: based on the theory of affine spheres, we
propose a new definition of area distance on the outer part of a
convex curve and new algorithms for computing these distances.
This link is also interesting from the point of view of the theory
of improper affine spheres, since it provides a geometrical
interpretation of them.

Let us make this connection more precise: we first remark that the
area distance $f$ to a convex plane curve satisfies the
Monge-Amp\`{e}re differential equation $det(D^2(f))=-1$ with
boundary conditions $f=0$ and $\nabla(f)=0$ (see \cite{Alvim05}).
We show in this paper that the graph of $f$ is an indefinite
improper affine sphere, with strictly positive Pick invariant. On
the other hand, we show also that, at least locally, any
indefinite improper affine sphere is an area distance.

For the outer part of the convex curve, we propose to define an
area based distance by the Monge-Amp\`{e}re differential equation
$det(D^2(f))=+1$ with boundary conditions $f=0$ and $\nabla(f)=0$.
Then the graph of $f$ is a definite improper affine sphere, also
with strictly positive Pick invariant. We consider in this paper
the case of an initial analytic curve. In this case, it is
possible to find explicitly a solution to the Monge-Amp\`{e}re
equation, and also to describe its geometrical properties.

Consider now a convex polygon as a discretization of the convex
curve. By following the asymptotic lines, we propose a fast
evolution algorithm that computes exactly the area distance. We
also show that this exact area distance defines a discrete
indefinite improper affine sphere, as defined in
\cite{Matsuura03}. We remark that the method proposed in this
paper differs completely from that of \cite{Alvim07}, since the
latter considers curves defined in implicit form. For the outer
part of the polygon, we also propose a fast evolution algorithm
that computes a new distance. This graph of this new distance has
the remarkable property of being a discrete definite improper
affine sphere, as defined in \cite{Matsuura03}.

In this context, there is a natural duality between points in the
inner and the outer part of the convex curve. An interesting fact
is that, although this duality is not area preserving, it
preserves the measure $J^{1/3}dxdy$, where $J$ denotes the Pick
invariant of the corresponding graph. From the discrete point of
view, it is interesting to observe that a mesh with planar crosses
in the inner part of the curve changes smoothly along the curve to
a mesh with planar quadrilaterals in the outer part.

This paper is organized as follows: In section 2, we review the
basic concepts related to improper affine spheres, both smooth and
discrete. In section 3, we review the definition of area distance
in the inner part of a curve and show its strong link with
indefinite improper affine spheres. In section 4, we propose the
new definition of area distance in the outer part of the curve,
show that its graph is an definite improper affine sphere and
describe the duality between the inner and the outer area
distances.

\smallskip\noindent{\bf Notation.} For three vectors $X,Y$ and $Z$ in the space, denote by
$[X,Y,Z]$ the determinant of the $3\times 3$ matrix whose columns
are the vectors $X,Y$ and $Z$. For two vectors $X$ and $Y$ in the
plane, denote by $[X,Y]$ the determinant of the $2\times 2$ matrix
whose columns are the vectors $X$ and $Y$. Also, denote by $X^t$
the transpose of the matrix $X$ and by $R$ the ninety degrees
rotation in the anti-clockwise direction. Observe that
$[X,Y]=-X^t\cdot R\cdot Y$.

\section{Improper affine spheres}

\subsection{Basics}

For a surface $S$ parameterized  by
$q(u,v)=(x(u,v),y(u,v),z(u,v))$ let $L=[q_u,q_v,q_{uu}]$,
$M=[q_u,q_v,q_{uv}]$ and $N=[q_u,q_v,q_{vv}]$. We say that $S$ is
non-degenerate if $LN-M^2\neq 0$. For a non-degenerate surface,
the Blaschke metric is given by
$$
\phi=\frac{Ldu^2+2Mdudv+Ndv^2}{|LN-M^2|^{1/4}}.
$$
It is definite or indefinite according to $LN-M^2$ being positive
or negative. In the definite case, the affine normal is defined as
$\xi=\frac{1}{2}\Delta(q)$, while in the indefinite case
$\xi=-\frac{1}{2}\Delta(q)$, where $\Delta(q)$ denotes the
laplacian of each coordinate with respect to the Blaschke metric.

An improper affine sphere $S$ is a surface whose affine normals at
all points are parallel. We shall assume that the affine normals
are parallel to the $z$-axis. Under this hypothesis $S$ is locally
the graph of a function $z=f(x,y)$. The following proposition is
well-known (see \cite{Simon93}).

\begin{proposition} $S$ is an indefinite improper affine sphere if and only if
$det(D^2(f))=-c$, and a definite improper affine sphere if and
only if $det(D^2(f))=c$, where $c$ is a positive constant. In both
cases, the affine normal is constant and equal to $(0,0,c)$.
\end{proposition}

In the case $S$ is the graph of $f$, we shall write
$q(u,v)=(p(u,v),f(p(u,v)))$. The coefficients of the Blaschke
metric can be calculated by the following lemma, whose proof is a
straightforward calculation:

\begin{lemma}\label{lema1}
$$
 \left\{
\begin{array}{l}
L=[p_u,p_v] D^2(f)(p_u,p_u)\\
M=[p_u,p_v] D^2(f)(p_u,p_v)\\
N=[p_u,p_v] D^2(f)(p_v,p_v)
\end{array}\ .
\right.
$$
\end{lemma}

\subsection{ Asymptotic and isothermal directions}

For an indefinite Blaschke metric, one can find parameters $(u,v)$
such that $L=N=0$. These parameters are called {\it asymptotic
parameters} and the corresponding tangent vectors are called {\it
asymptotic directions}. Lines whose tangent vectors are asymptotic
directions are called {\it asymptotic lines}. Lemma \ref{lema1}
shows that the projections of the asymptotic directions in the
$(x,y)$-plane vanish the quadratic form $D^2(f)$, and we shall
also call them asymptotic directions.

 Using asymptotic parameters, the structure equations become
\[
\left\{
\begin{array}{l}
q_{uu}=\frac{\omega_u}{\omega} q_u+\frac{a}{\omega} q_v \\
q_{vv}=\frac{b}{\omega} q_u+\frac{\omega_v}{\omega} q_v \\
q_{uv}=-\omega\xi
\end{array}
\right.
\]
where $a_v=b_u=0$ and $\omega=[q_u,q_v,\xi]$. By a good choice of
the asymptotic parameters, we can make $a$ and $b$ constants. The
Pick invariant is given by $j=\frac{ab}{\omega^3}$ (see
\cite{Buchin83} and \cite{Matsuura03}). The equations for the
planar component are
\[
\left\{
\begin{array}{l}
p_{uu}=\frac{\omega_u}{\omega} p_u+\frac{a}{\omega} p_v \\
p_{vv}=\frac{b}{\omega} p_u+\frac{\omega_v}{\omega} p_v \\
p_{uv}=0,
\end{array}
\right.
\]
with $\omega=[p_u,p_v]=-f_{uv}$. From these equations one obtains
$[p_{u},p_{uu}]=a$ and $[p_{v},p_{vv}]=-b$. So, if the Pick
invariant does not vanish, the planar asymptotic lines are convex.

For a definite Blaschke metric, one can consider also asymptotic
parameters, but they are complex (\cite{Buchin83}).  Consider
complex parameters $z=s+it$ and $\zb=s-it$ such that
$L(z,\zb)=N(z,\zb)=0$ and $iM(z,\zb)=\Omega^2$. Then, in terms of
$s$ and $t$, $L(s,t)=N(s,t)=4\Omega^2$ and $M(s,t)=0$
(\cite{Bobenko199}). We call such coordinates {\it isothermal}.

 In order to
differentiate from the indefinite case, we shall use capital
letters $P,Q,F,\Omega,A,B$ and $J$ to describe the structure
equations, as follows:
\[
\left\{
\begin{array}{l}
Q_{zz}=\frac{\Omega_z}{\Omega} Q_z-\frac{A}{\Omega} Q_{\zb} \\
Q_{\zb\zb}=-\frac{B}{\Omega} Q_z+\frac{\Omega_{\zb}}{\Omega} Q_{\zb} \\
Q_{z\zb}=-\Omega\xi
\end{array}
\right.
\]
where $\Omega=-i[Q_z,Q_{\zb},\xi]$, $A=i[Q_z,Q_{zz},\xi]$ and
$B=-i[Q_{\zb},Q_{\zb\zb},\xi]$. By a good choice of the isothermal
parameters, we can make $A$ and $B$ constants. The Pick invariant
is given by $J=\frac{AB}{\Omega^3}$ (see \cite{Bobenko199} and
\cite{Matsuura03}). The equations for the planar components are
\[
\left\{
\begin{array}{l}
P_{zz}=\frac{\Omega_z}{\Omega} P_z-\frac{A}{\Omega} P_{\zb} \\
P_{\zb\zb}=-\frac{B}{\Omega} P_z+\frac{\Omega_{\zb}}{\Omega} P_{\zb} \\
P_{z\zb}=0, \ \ F_{z\zb}=-\Omega
\end{array}
\right.
\]
with $\Omega=-i[P_z,P_{\zb}]$. Also $A=i[P_{z},P_{zz}]$ and
$B=-i[P_{\zb},P_{\zb\zb}]$.

\subsection{Discrete improper affine spheres.}

In \cite{Bobenko199}, definitions of discrete proper affine
spheres are proposed, both in the indefinite and in the definite
case. In \cite{Matsuura03}, these definitions are generalized to
discrete improper affine spheres, indefinite and definite. We
describe now these latter definitions, with a slight modification
in the definite case. Denote by $Z$ the set of integers.

\begin{definition} \label{IndefiniteDiscrete} A map $q:Z^2\to R^3$ is a discrete
indefinite improper affine sphere if it has the following
properties:
\begin{enumerate}
\item For any $(i,j)\in Z^2$, the points
$q(i,j),q(i+1,j),q(i-1,j),q(i,j-1),q(i,j+1)$ are co-planar.
\item There exists a direction $\xi$ in $R^3$ such that, for any $(i,j)\in Z^2$,
the vector $q(i,j)+q(i+1,j+1)-q(i+1,j)-q(i,j+1)$ is parallel to
$\xi$.
\end{enumerate}
\end{definition}

\begin{definition}\label{DefiniteDiscrete} A map $Q:Z^2\to R^3$ is a discrete definite
improper affine sphere if it has the following properties:
\begin{enumerate}
\item For any $(i,j)\in Z^2$, the points
$Q(i,j),Q(i+1,j),Q(i,j+1),Q(i+1,j+1)$ are co-planar.
\item There exists a direction $\xi$ in $R^3$ such that, for any $(i,j)\in Z^2$,
the vector $Q(i,j+1)+Q(i,j-1)+Q(i-1,j)+Q(i+1,j)-4Q(i,j)$ is
parallel to $\xi$.
\end{enumerate}
\end{definition}

\section {Inner Area distances}

In this section we review some properties of area distances and
show the connection between area distances and affine spheres. We
show that the graph of an area distance is an indefinite improper
affine sphere and that, at least locally, any indefinite improper
affine sphere is the graph of an area distance of some convex
plane curve.

Moreover, we show that the area distance of polygons define
discrete indefinite improper affine spheres and how this fact can
be applied to construct a very fast algorithm for computing area
distances.

\subsection{Definition and properties of the area distance}

\begin{figure}[htb]
\centering \fsep\subfigure[Domain for a curve with 2 endpoints. ]
{\includegraphics[width=.25 \linewidth,clip
=false]{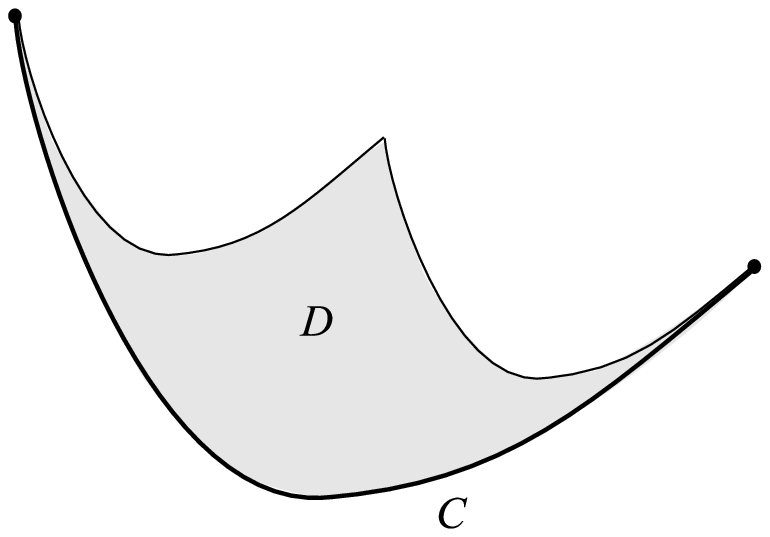}} \fsep\fsep\subfigure[ Domain for a curve
with 1 endpoint. ] {
\includegraphics[width=.25\linewidth,clip
=false]{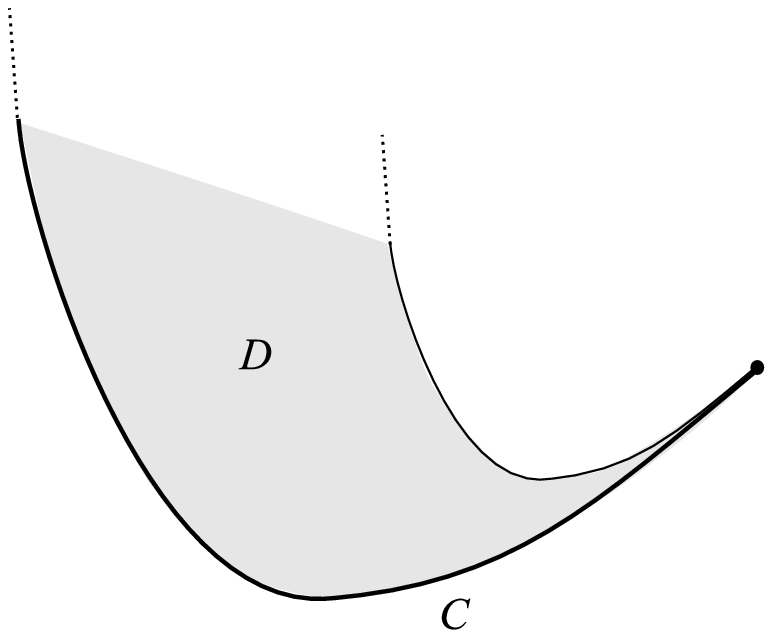}}\fsep\fsep \subfigure[Domain for a curve
without endpoints. ] {\includegraphics[width=.25 \linewidth,clip
=false]{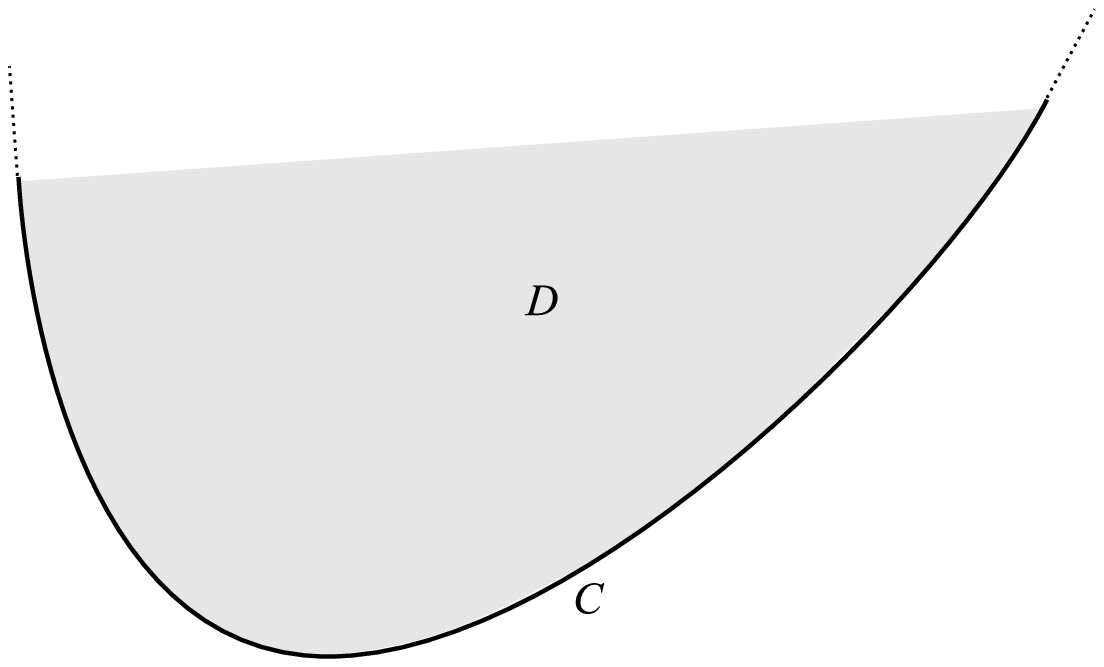}}\fsep\caption{ Domains for the inner area
distance function.} \label{dominio}
\end{figure}

Consider a smooth convex curve $C$ in the plane without parallel
tangent lines. $C$ can have $2$, $1$ or $0$ endpoints. Denote by
$D$ the plane region whose boundary is $C$ and the curve(s)
obtained from $C$ by a similarity of ratio $\frac{1}{2}$ based at
each endpoint of $C$ (see figure \ref{dominio}).

 A chord is a line segment connecting $2$ points of $C$.
For a point $p\in D$, consider chords $l$ passing through $p$
that, together with $C$, bound regions $D_l$. Denote by $l(p)$ the
chord such that the area of the corresponding region $D_{l(p)}$ is
minimum. The area distance function $f(p)$ is then defined as half
of the area of $D_{l(p)}$. Sometimes we shall call this function
{\it inner area distance}, since in section \ref{exterior} we
shall define another area distance.

\begin{figure}[htb]
 \centering
 \includegraphics[width=0.35\linewidth]{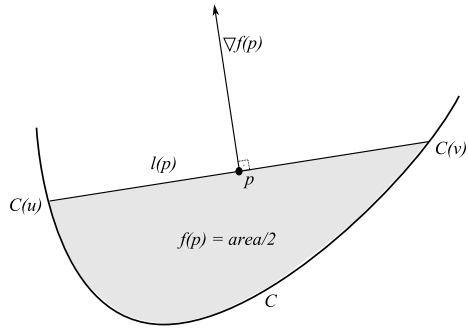}
 \caption{Gradient of the area distance.}
\label{GradientAreaDistance}
\end{figure}

Any $p\in D$ is the mean point of the extremities of the chord
$l(p)$. This important property was proved first in
\cite{Moisan98} (see also \cite{Giblin04}). Another important
property is that $\nabla(f)(p)$ is orthogonal to the chord $l(p)$,
with half of the length of it (see figure
\ref{GradientAreaDistance}). The third important property is that
$det(D^2(f)(p))=-1$, for any $p\in D$. This last property was
first proved in \cite{Alvim05}. The next lemma and proposition
describe the $2$ latter properties with more details.

\begin{lemma}
Denote by $C(u(p))$ and $C(v(p))$ the extremities of the chord
$l(p)$. Considering $u$ and $v$ as functions of $p$, we have $
\nabla(u)=-\frac{2R C'(v)}{[C'(u),C'(v)]} $ and $
\nabla(v)=\frac{2R C'(u)}{[C'(u),C'(v)]} $.
\end{lemma}

\proof{ Since $2p=C(u(p))+C(v(p))$, $
2I=C'(u)\nabla(u)^t+C'(v)\nabla(v)^t $ where $I$ is the identity
$2\times 2$ matrix. Multiplying by $ C'(v)^t R$, one obtains the
first formula. The second one is analogous.}

\begin{proposition}\label{GradienteHessiana} The area distance $f$ of a convex arc $C$
satisfies the following formulas:

\begin{enumerate}

\item   The gradient of $f$ is given by $\nabla(f)=R(\frac{1}{2}(C(v)-C(u)))$.

\item \label{P2} $C'(u)$ and $C'(v)$ satisfy the equations
$D^2(f)(C'(u),C'(u))=0$, $D^2(f)(C'(v),C'(v))=0$ and
$D^2(f)(C'(u),C'(v))=[C'(u),C'(v)]$.

\item   $det(D^2(f))=-1$.

\end{enumerate}

\end{proposition}
\proof{

\begin{enumerate}
\item
 By Green's theorem, one can write
$4f(p)=\int_{u}^{v}[C(s)-p,C'(s)]ds$. Using the rules for
differentiating integrals, one obtains
$$
4\nabla(f)(p)=[C(v)-p,C'(v)]\nabla(v)-[C(u)-p,
C'(u)]\nabla(u)+\int_{u}^{v}R C'(s)ds,
$$
$$
4\nabla(f)(p)=R (\frac{[C(v)-C(u),C'(v)]
}{[C'(u),C'(v)]}C'(u)+\frac{[C(v)-C(u), C'(u)]
}{[C'(v),C'(u)]}C'(v))+R(C(v)-C(u)).
$$
Since, for any $w$, one can write
$w=\frac{[C'(v),w]}{[C'(v),C'(u)]}C'(u)+\frac{[C'(u),w]}{[C'(u),C'(v)]}C'(v)$,
we conclude that
$$
\nabla(f)(p)=R(\frac{1}{2}(C(v)-C(u))).
$$
\item
 Differentiating  one obtains
$$
2D^2(f)=R(C'(u)\nabla(u)^t-C'(v)\nabla(v)^t),
$$
and so
$$
D^2(f)=\frac{1}{[C'(u),C'(v)]}\left(RC'(u)C'(v)^t R+R C'(v)C'(u)^t
R\right).
$$
One concludes that
$$
C'(u)^t D^2(f) C'(u)=C'(v)^t D^2(f) C'(v)=0\
$$
and
$$
C'(u)^t D^2(f) C'(v)=C'(v)^t D^2(f) C'(u)=[C'(u),C'(v)].
$$
\item Since we have assumed that there are no parallel tangents, property \ref{P2} guarantees that $D^2(f)$ is non-degenerate. Now
$C'(u)^tD^2(f)C'(u)=0$ implies that $D^2(f)C'(u)=\lambda RC'(u)$,
for some $\lambda\neq 0$. And since
$C'(v)^tD^2(f)C'(u)=[C'(u),C'(v)]$, $\lambda=1$. So
$D^2(f)C'(u)=RC'(u)$ and similarly $D^2(f)C'(v)=-RC'(v)$. Hence
$[D^2(f)C'(u),D^2(f)C'(v)]=[RC'(u),-RC'(v)]=-[C'(u),C'(v)]$, which
implies that $det(D^2(f))=-1$.
\end{enumerate}}

\subsection{Area distances as improper affine spheres}

\subsubsection{Consequences of proposition
\ref{GradienteHessiana}}

Proposition \ref{GradienteHessiana} implies that  the graph of the
area distance $f$ to a smooth convex arc $C$ is an indefinite
improper affine sphere. Also, the parameterization
$$
q(u,v)=(\frac{1}{2}(C(u)+C(v)),f(\frac{1}{2}(C(u)+C(v)))).
$$
is asymptotic. Hence $C'(u)$ and $C'(v)$ are asymptotic directions
and the asymptotic lines are obtained from $C$ by similarities of
ratio $1/2$. Direct calculations show that
$\omega(u,v)=\frac{1}{4}[C'(v),C'(u)]$. Finally, $f(p)$ represents
the area of the region bounded by the two asymptotic lines that
start at $p$ and the curve $C$ itself (see figure
\ref{projections}).

\smallskip
 If the parameterization of the curve $C$ is by affine
arc length, then $a=b=-\frac{1}{4}$. In this case the Pick
invariant is $j=4[C'(v),C'(u)]^{-3}$. Thus it vanishes only at
points $p$ such that at least one of the endpoints of minimal
chord $l(p)$ belongs to a line segment of the original curve $C$.
In particular, if the curve is strictly convex, the Pick invariant
never vanishes.

\subsubsection{Some explicit formulas for calculating $f$.}

By the first item of proposition \ref{GradienteHessiana}, formulas
$f_{u}=f_{x}x_{u}+f_{y}y_{u}$ and $f_{v}=f_{x}x_{v}+f_{y}y_{v}$
imply that

\begin{eqnarray*}
f_u\left( u,v\right)  &=&\frac{1}{4}\left\vert
\begin{array}{cc}
x\left( u\right) -x\left( v\right)  & x^{\prime }\left( u\right)  \\
y\left( u\right) -y\left( v\right)  & y^{\prime }\left( u\right)
\end{array}%
\right\vert  \\
f_v\left( u,v\right)  &=&\frac{1}{4}\left\vert
\begin{array}{cc}
x\left( u\right) -x\left( v\right)  & x^{\prime }\left( v\right)  \\
y\left( u\right) -y\left( v\right)  & y^{\prime }\left( v\right)
\end{array}%
\right\vert
\end{eqnarray*}%

It is also interesting to consider coordinates $(s,t)$ defined by
$2s=u+v$ and $2t=u-v$. In these coordinates, since $g_s=g_u+g_v$
and $g_t=g_u-g_v$, we have

\begin{eqnarray}
f_s\left( s,t\right)  &=&\frac{1}{4}\left\vert
\begin{array}{cc}\label{fs}
x\left( s+t\right) -x\left( s-t\right)  & x^{\prime }\left( s+t\right)+ x^{\prime }\left( s-t\right) \\
y\left( s+t\right) -y\left( s-t\right)  & y^{\prime }\left(
s+t\right)+y^{\prime }\left( s-t\right)
\end{array}%
\right\vert  \\
f_t\left( s,t\right)  &=&\frac{1}{4}\left\vert
\begin{array}{cc}\label{ft}
x\left( s+t\right) -x\left( s-t\right)  & x^{\prime }\left( s+t\right) -x^{\prime }\left( s-t\right) \\
y\left( s+t\right) -y\left( s-t\right)  & y^{\prime }\left(
s+t\right)-y^{\prime }\left( s-t\right)
\end{array}%
\right\vert
\end{eqnarray}%

\subsection{Examples}

In this subsection we give explicit examples of area based
distance function of convex smooth curves.

\begin{example}\label{parabola}
 Consider the parabola parameterized by
$C\left( r\right)=\left( r,\frac{r^{2}}{2}\right)$. Integrating
$(f_u,f_v)=\frac{1}{8}((u-v)^2,-(u-v)^2)$ one obtains
$f(u,v)=\frac{1}{24}(u-v)^3$. Thus a parameterization of the
affine sphere in asymptotic coordinates is given by
$$
q(u,v)=(\frac{u+v}{2},\frac{u^2+v^2}{4},\frac{1}{24}(u-v)^3)
$$
with $u>v$. And $\omega(u,v)=\frac{1}{4}(u-v)$. In $(s,t)$
coordinates,
$$
q(s,t)=(s,\frac{s^2+t^2}{2},\frac{t^3}{3})
$$
and $\omega(s,t)=\frac{1}{2}t$.

\end{example}

\begin{example}\label{circle}
 Consider the circle parameterized by $C(r)=(\cos(r),\sin(r))$. Although this curve admits parallel tangents,
the above calculations work well, except at the center of the
circle. Integrating $(f_u,f_v)=\frac{1}{4}(-\cos \left( u-v\right)
+1,\cos \left( u-v\right) -1)$ one obtains
$f(u,v)=\frac{1}{4}(u-v-\sin \left( u-v\right))$. Thus a
parameterization in asymptotic coordinates is given by
$$
q(u,v)=(\frac{1}{2}(\cos(u)+\cos(v)),\frac{1}{2}(\sin(u)+\sin(v)),\frac{1}{4}(u-v-\sin(u-v)))
$$
with $u>v$. And $\omega(u,v)=\frac{1}{4}\sin(u-v)$. In $(s,t)$
coordinates,
$$
q(s,t)=(\cos(s)\cos(t),\sin(s)\cos(t),\frac{1}{4}(2t-\sin(2t)))
$$
and $\omega(s,t)=\frac{1}{4}sin(2t)$.

\end{example}
\begin{example}\label{hyperbola}
Consider the hyperbola parameterized by $C(r)=(\exp(r),\exp(-r))$.
Integrating
$(f_u,f_v)=\frac{1}{4}(-2+\exp(u-v)+\exp(v-u),2-\exp(u-v)-\exp(v-u))$
one obtains
$f(u,v)=\frac{1}{2}(v-u)+\frac{1}{4}(\exp(u-v)-\exp(v-u))$. Thus
the asymptotic parameterization is
$$
q(u,v)=\frac{1}{2}(\exp(u)+\exp(v),\exp(-u)+\exp(-v),(v-u)+\frac{1}{2}(\exp{(u-v)}-\exp{(v-u)})).
$$
with $u>v$. And $\omega(u,v)=\frac{1}{4}(\exp(u-v)-\exp(v-u))$. In
$(s,t)$ coordinates,
$$
q(s,t)=\frac{1}{2}(\exp(s+t)+\exp(s-t),\exp(-(s+t))+\exp(t-s),-2t+\frac{1}{2}(\exp{(2t)}-\exp{(-2t)})).
$$
and $\omega(s,t)=\frac{1}{4}(\exp(2t)-\exp(-2t))$.

\end{example}

\begin{example}\label{cubic}
Consider the cubic parameterized by $C(r)=(r,r^3)$. Although this
is not a convex curve, the above calculations can be done.
Integrating the vector field
$(f_u,f_v)=\frac{1}{4}(2u^{3}-3u^{2}v+v^{3},u^{3}-3uv^{2}+2v^{3}
)$ one obtains $f(u,v)=\frac{1}{8}\left( u-v\right) ^{3}\left(
u+v\right)$. Thus the asymptotic parameterization is
$$
q(u,v)=\frac{1}{2}((u+v),u^3+v^3,\frac{1}{4}\left( u-v\right)
^{3}\left( u+v\right)).
$$
with $u>v$. And $\omega(u,v)=\frac{3}{4}(u^2-v^2)$. In $(s,t)$
coordinates,
$$
q(s,t)=(s,s(3t^2+s^2),2t^3s)
$$
and $\omega(s,t)=3st$.

\end{example}

\begin{example}\label{r4}
Consider the curve parameterized by $C(r)=(r,r^4)$. Integrating
the vector field
$(f_u,f_v)=\frac{1}{4}(3u^4-4u^3v+v^4,4v^3u-3v^4-u^4)$ one obtains
$f(u,v)=\frac{1}{4}(\frac{3u^5}{5}-\frac{3v^5}{5}-u^4v+uv^4)$.
Thus the asymptotic parameterization is
$$
q(u,v)=\frac{1}{2}(u+v,u^4+v^4,\frac{1}{2}(\frac{3u^5}{5}-\frac{3v^5}{5}-u^4v+uv^4)
).
$$
with $u>v$. And $\omega(u,v)=u^3-v^3$. In $(s,t)$ coordinates,
$$
q(s,t)=(s,t^4+s^4+6s^2t^2,4t^3s^2+\frac{4t^5}{5})
$$
and $\omega(s,t)=2t^3+6ts^2$.

\end{example}

\subsection{Local characterization of indefinite improper affine
spheres}\label{generality}

In this subsection, we show that locally, and up to a constant,
any indefinite improper affine sphere with non-zero Pick invariant
is the graph of the area distance of a smooth convex plane curve
$C$. More precisely, we have the following theorem:

\begin{theorem}\label{LocalCharacterization}
Let $U$ be an open domain in the $(u,v)$-plane whose closure
$\overline{U}$ is contained in the domain of the asymptotic
parameterization of an indefinite improper affine sphere $S$.
Assume that, restricted to $U$, $S$ is the graph of a function
$f(p)$. Then there exists a convex curve $C$ in the plane and a
constant $K$ such that $f+K$ is the area distance of $C$.

\end{theorem}

\subsubsection{Some properties of indefinite improper affine spheres}

Consider an indefinite improper affine sphere with strictly
positive Pick invariant $j$. Assume, w.l.o.g., that the affine
normal is $\xi=(0,0,1)$ and consider that $S$ is the graph of a
function $f$.

\begin{figure}[htb]
\centering
\includegraphics[width=0.35\linewidth]{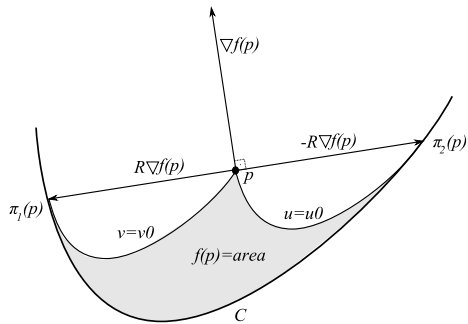}
\caption{The projections $\pi_1$ and $\pi_2$.} \label{projections}
\end{figure}

\begin{lemma}\label{properties} The following properties hold (see figure \ref{projections}):

\begin{enumerate}
\item We have that
$D^2(f)(p_u)=Rp_u $ and $D^2(f)(p_v)=-Rp_v$.
\item Let $\pi_1(p)=p+R\nabla(f)(p)$ and $\pi_2(p)=p-R\nabla(f)(p)$. Then
$\pi_1(p)$ is constant along an asymptotic line $v=v_0$ and
$\pi_2(p)$ is constant along an asymptotic line $u=u_0$.
\end{enumerate}
\end{lemma}

 \proof{
\begin{enumerate}
\item
Since $D^2(f)$ is non-degenerate and $D^2(f)(p_u,p_u)=0$,
$D^2(f)p_u=\lambda Rp_u$, for some $\lambda\neq 0$. And since
$D^2(f)(p_u,p_v)=\omega(u,v)=[p_u,p_v]$, $\lambda=1$. A similar
reasoning shows that $D^2(f)(p_v)=-Rp_v$.
\item
Just observe that $p_u+RD^2(f)p_u=0$ and $p_v-RD^2(f)p_v=0$.
\end{enumerate}
}

\subsubsection{Proof of theorem \ref{LocalCharacterization} }

Let  $\pi_1$ and $\pi_2$ be the projections defined in lemma
\ref{properties}. We have that $C_1=\pi_1(\overline{U})$ and
$C_2=\pi_2(\overline{U})$ are compact arcs. It is not difficult to
obtain a smooth convex arc $C_3$ such that the concatenation $C$
of $C_1$, $C_3$ and $C_2$ is smooth nd convex. Denote by $g$ the
area distance function associated to $C$.

One can easily see that $p=\frac{1}{2}(\pi_1(p)+\pi_2(p))$ and
$\nabla(f)(p)=\frac{1}{2}R(\pi_2(p)-\pi_1(p))$, for any $p\in
\overline{U}$. So $\nabla(g)(p)=\nabla(f)(p)$, for any $p\in
\overline{U}$. This implies that $f-g$ is constant, which proves
the theorem.

\subsection{Area distances to polygons}

\subsubsection{Asymptotic grids 2d and 3d}

\begin{figure}[htb]
\centering \fsep\subfigure[Points and parallelograms. ]
{\includegraphics[width=.25 \linewidth,clip
=false]{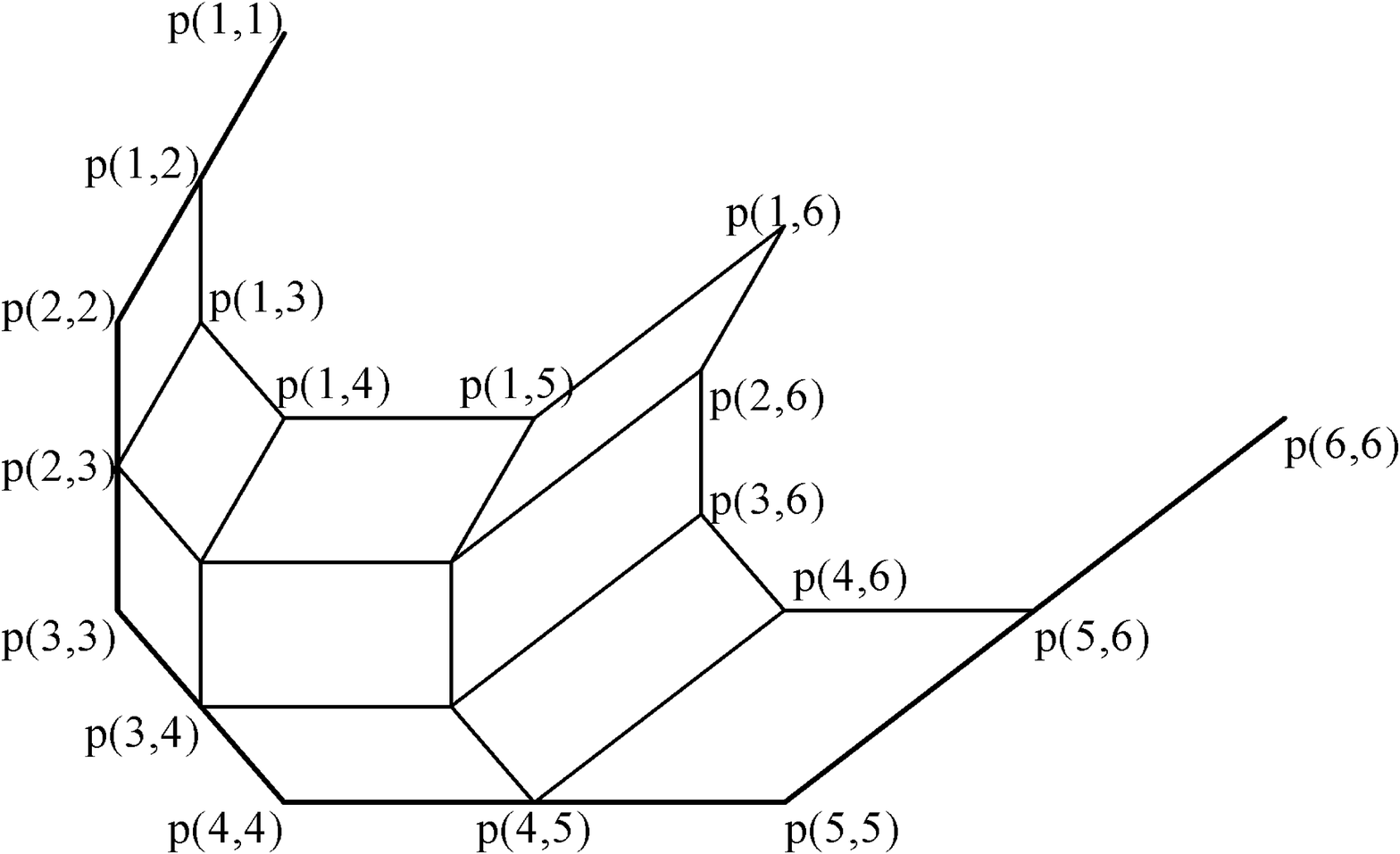}} \fsep\fsep\subfigure[ Areas defined by the
parallelograms. ] {
\includegraphics[width=.25\linewidth,clip
=false]{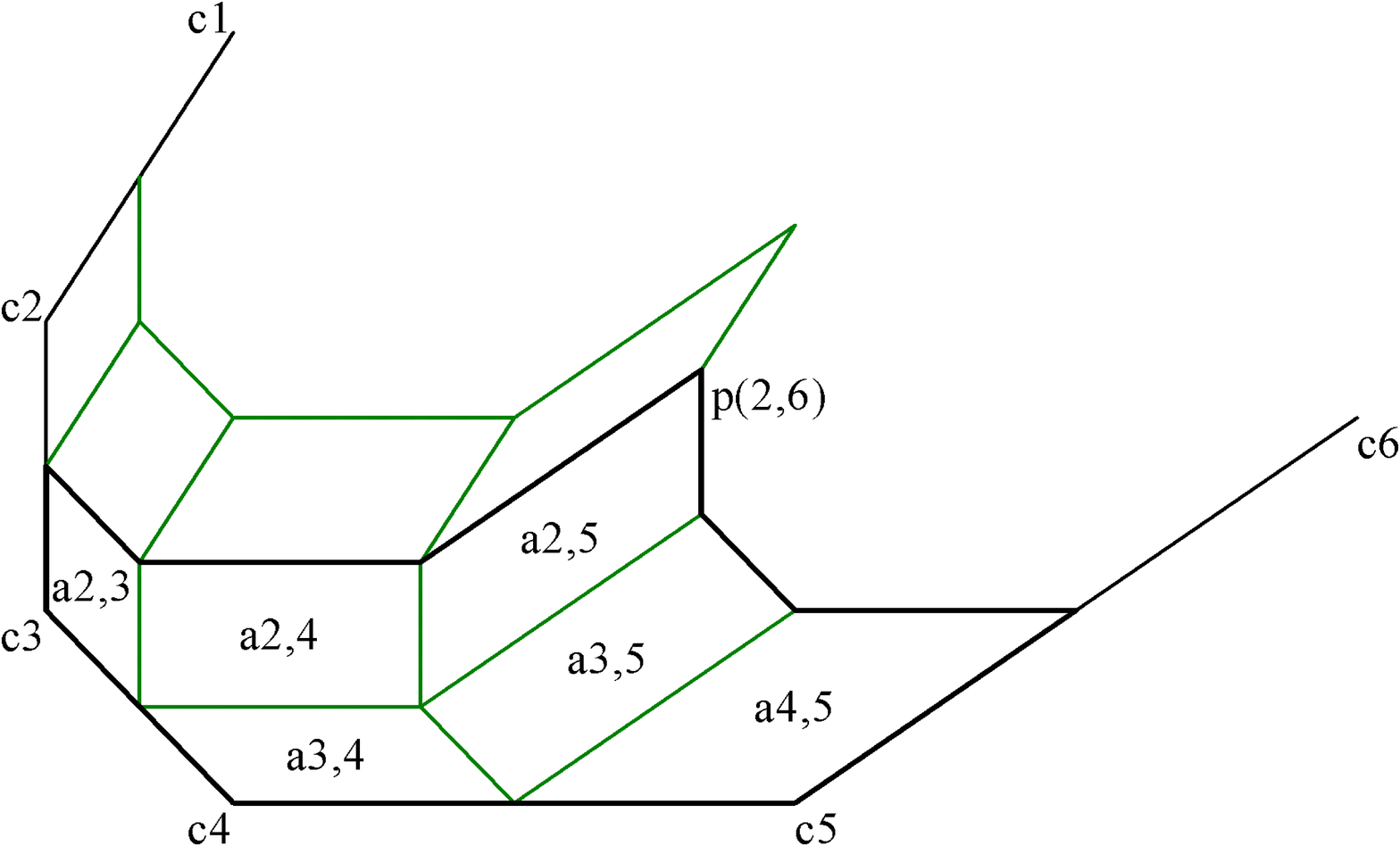}}\fsep\fsep\fsep \caption{ Planar mesh of a
discrete inner distance.} \label{areas}
\end{figure}

Let $C$ be a convex polygon with vertices $c_i$, $0\leq i\leq
N-1$. Denote its sides by the vectors $2L_{i}=c_{i+1}-c_{i}$ and
assume that $[L_i,L_j]>0$, for any $j>i$.
Assuming $j\geq i$, define the following grid $%
p\left( i,j\right) $ on the plane by
\[
p\left( i,j\right) =\frac{c_{i}+c_{j}}{2}.
\]%
Note that%
\[
p\left( i+1,j\right) -p\left( i,j\right) =\frac{c_{i+1}-c_{i}}{2}%
=2L_{i}=p\left( i+1,j+1\right) -p\left( i,j+1\right),
\]%
so the grid is formed by parallelograms whose areas will be denoted by%
\[
a_{ij}=\left[ p\left( i+1,j\right) -p\left( i,j\right) ,p\left(
i,j+1\right) -p\left( i,j\right) \right] =\left[
L_{i},L_{j}\right].
\]
(see figure \ref{areas}). Define also
\[
f\left( i,j\right) =\sum_{i\leq k<l\leq j-1}a_{k,l}=\sum_{i\leq
k<l\leq j-1} \left[ L_{k},L_{l}\right].
\]
Note that $f\left( i,j\right) =0$ if $j=i$ or $j=i+1$.



\begin{proposition}

The map $q\left( i,j\right) =\left( p\left( i,j\right) ,f\left(
i,j\right) \right) $ is a discrete indefinite improper affine
sphere.

\end{proposition}

\begin{figure}[htb]
 \centering
 \includegraphics[width=0.35\linewidth]{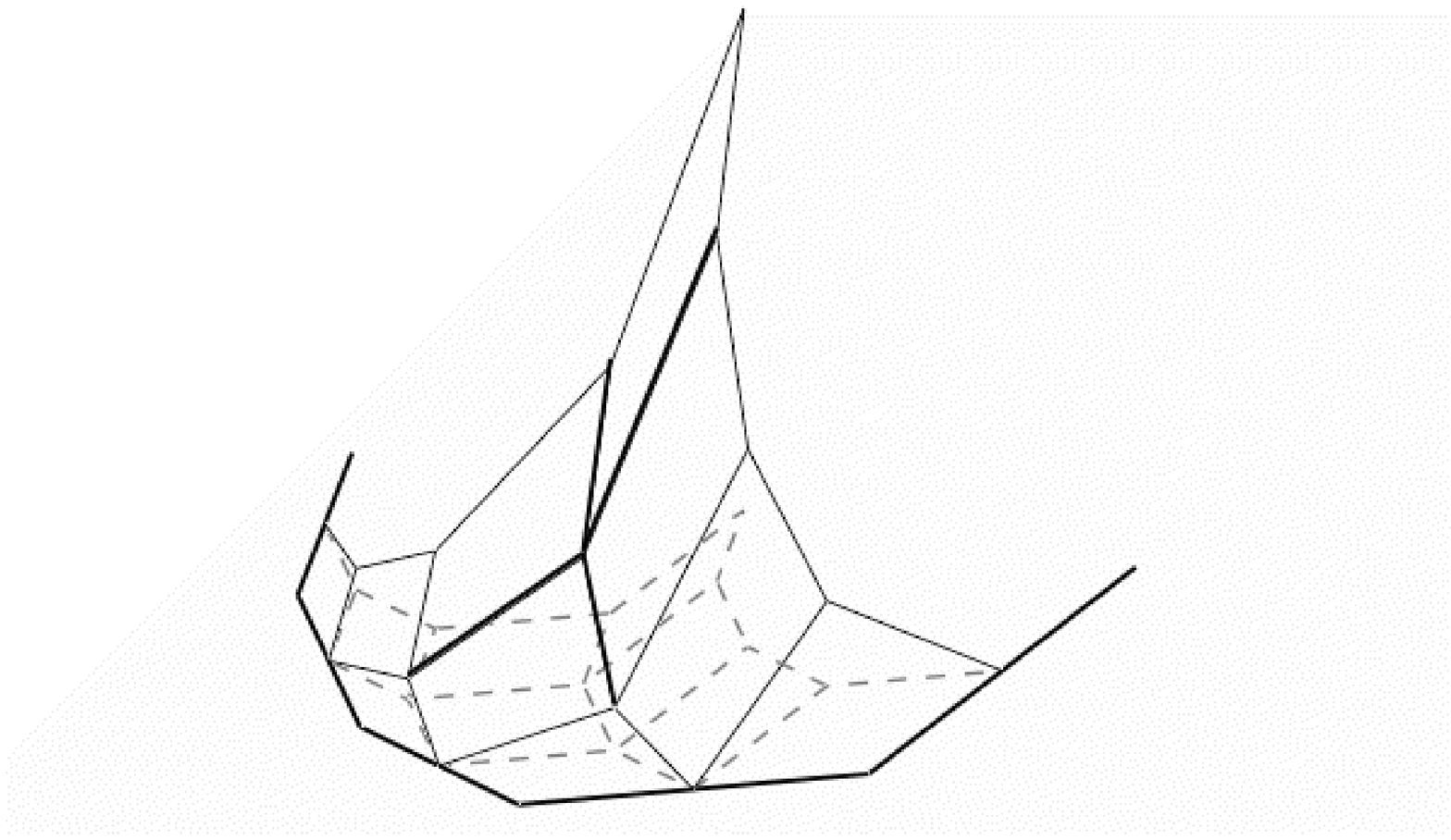}
 \caption{The four line segments in bold are co-planar.}
\label{cross}
\end{figure}

\proof{

\begin{enumerate}
\item Let%
\[
L=\sum_{k=i}^{j-1}L_{k}
\]%
Then, since $\left[ L_{i},L_{i}\right] =\left[ L_{j},L_{j}\right] =0$:%
\begin{eqnarray*}
q\left( i,j\right) -q\left( i-1,j\right)  &=&\left(
L_{i-1},-\sum_{l=i}^{j-1} \left[ L_{i-1},L_{l}\right] \right)
=\left( L_{i-1},\left[ L,L_{i-1}\right]
\right)  \\
q\left( i+1,j\right) -q\left( i,j\right)  &=&\left(
L_{i},-\sum_{l=i+1}^{j-1} \left[ L_{i},L_{l}\right] \right)
=\left( L_{i},\left[ L,L_{i}\right]
\right)  \\
q\left( i,j+1\right) -q\left( i,j\right)  &=&\left(
L_{j},\sum_{k=i}^{j-1} \left[ L_{k},L_{j}\right] \right) =\left(
L_{j},\left[ L,L_{j}\right] \right)
\end{eqnarray*}%
so whatever linear dependence is satisfied by $L_{i-1}$, $L_{i}$
and $L_{j}$ will also be satisfied by the $z$ coordinates of the
three above vectors. This shows that $q\left( i+1,j\right) $ is in
the same plane as $q\left( i-1,j\right) $, $q\left( i,j-1\right) $
and $q\left( i,j\right) $. Similarly, one can show that $q\left(
i+1,j+1\right) $ is also in this plane (see figure \ref{cross}).

\item From the equations above, note that%
\begin{eqnarray*}
q\left( i+1,j+1\right) -q\left( i,j+1\right)  &=&\left(
L_{i},\left[
\sum_{k=i}^{j}L_{k},L_{i}\right] \right)  \\
q\left( i+1,j\right) -q\left( i,j\right)  &=&\left( L_{i},\left[
\sum_{k=i}^{j-1}L_{k},L_{i}\right] \right)
\end{eqnarray*}%
Subtracting,%
\[
q\left( i+1,j+1\right) -q\left( i,j+1\right) -q\left( i+1,j\right)
+q\left( i,j\right) =\left( 0,\left[ L_{j},L_{i}\right] \right)
=\left( 0,-a_{ij}\right)
\]%
and hence these vectors are all parallel to the $z$-axis.
\end{enumerate}
}

\subsubsection{Fast algorithm}

Now, if we define the level of a point to be $k=j-i$ we can
actually obtain our grid in levels starting from the polygon. At
level $0$, $p\left( i,i\right) =c_{i}$ and $f\left( i,i\right)
=0$. At level $1$, $p\left( i,i+1\right) =\frac{1}{2}\left(
c_{i}+c_{i+1}\right)$ and $f\left( i,i+1\right) =0$. At level $
\left( j-i\right) +1$,

\begin{eqnarray*}
p\left( i,j+1\right) =p\left( i,j\right) +p\left( i+1,j+1\right)
-p\left( i+1,j\right)  \\
f\left( i,j+1\right) =f\left( i,j\right) +f\left( i+1,j+1\right)
-f\left( i+1,j\right) +\left[ L_{i},L_{j}\right]
\end{eqnarray*}%

\begin{figure}[htb]
\centering \fsep\subfigure[ A discrete improper affine sphere. ] {
\includegraphics[width=.25
\linewidth,clip =false]{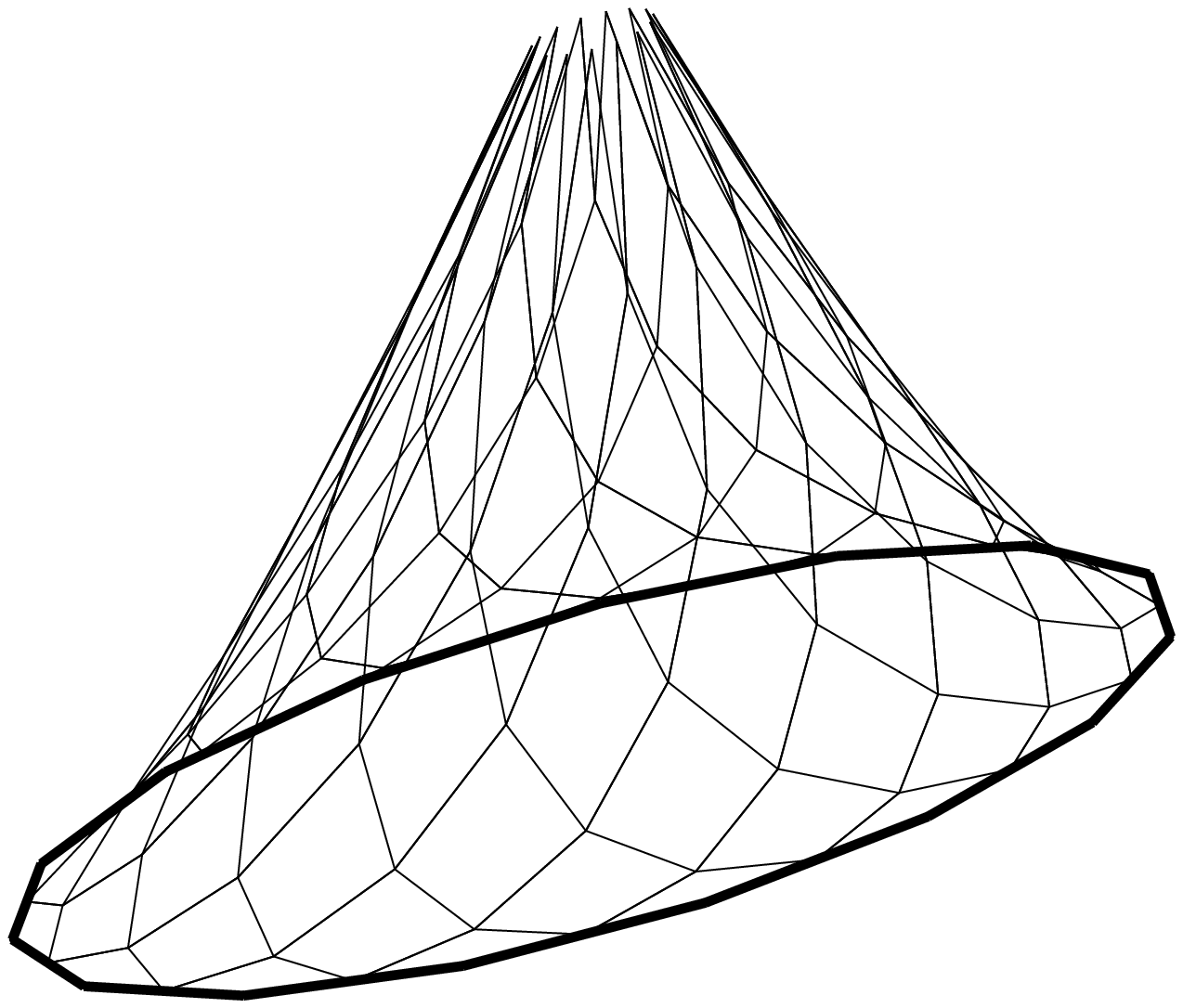}} \fsep\subfigure[ Planar
cross in bold. ] {
\includegraphics[width=.25\linewidth,clip
=false]{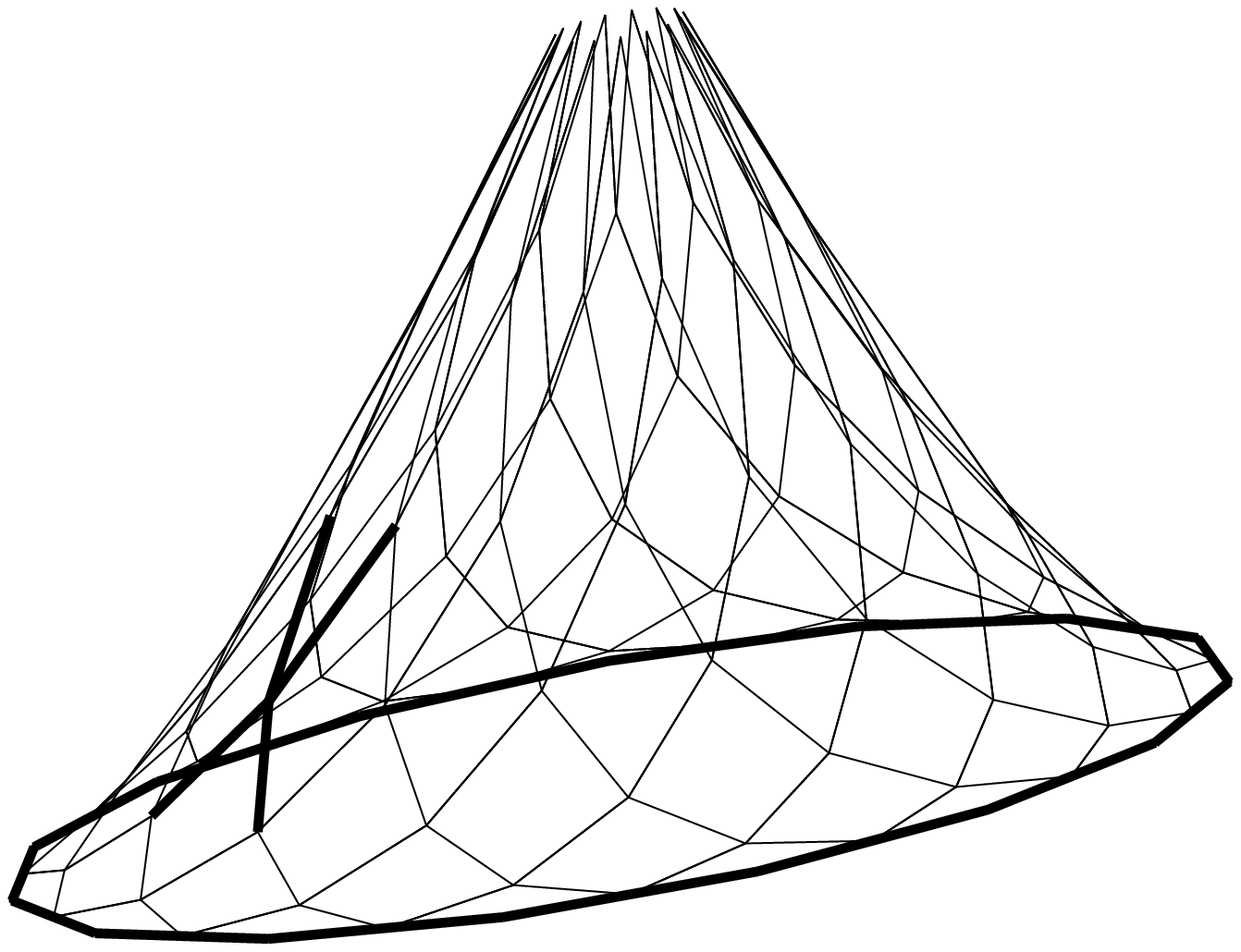}}\fsep \caption{ The inner
area distance of a polygon inscribed in a circle.}
\label{DiscreteInCircle3D}
\end{figure}

This gives us a fast algorithm to calculate all the grid points.
In each parallelogram, we can calculate the exact distance by a
bilinear interpolation. In figure \ref{DiscreteInCircle3D}, one
can see the result of this algorithm applied to a polygon
inscribed in an ellipse.

\section{Outer area distances }\label{exterior}

In this section, we consider the question of extending the area
distance to a neighborhood $E$ of $C$ in the outer part of $C$.
The idea is to solve the Monge-Amp\`{e}re differential equation
\[
\left\{
\begin{array}{l}
det(D^2F)(P)=+1, P\in E \\
\nabla F(P)=0, P\in C \\
F(P)=0, P\in C
\end{array}
\right.
\]
and define the area distance at $P$ by $F(P)$. It is clear that
the graph of $F$ defines an improper definite affine sphere.

In this section, we shall first describe a solution to this
problem in the case of an analytic curve. Then we prove some
properties of this solution, including a relation with the area
distance in $D$. Finally, we indicate how to obtain a discrete
definite improper affine sphere as an outer area distance of a
polygon.

\subsection{Definite improper affine spheres from analytic curves}

 Assume that the parameterization $C(r)=(x(r),y(r))$ is analytic
so that we can evaluate its coordinates for complex values of the
parameter, $z=s+it$, and $\zb=s-it$. Since $C$ is analytic and
real on the real line, the expression
$P(z,\zb)=\frac{1}{2}(C(z)+C(\zb))$ is real. It represents the
planar coordinates of the definite improper affine sphere defined
below.

\begin{lemma}
Consider the inner area function $f$ in coordinates $(s,t)$, and
let $F(s,t)=if(s,it)$. Then $F(s,t)$ is real and
$\nabla(F)=iR(\frac{C(z)-C(\zb)}{2})$.
\end{lemma}
\proof {From formulas (\ref{fs}) and (\ref{ft}), we can write

\begin{eqnarray}
F_s\left( s,t\right)  &=&\frac{i}{4}\left\vert
\begin{array}{cc}\label{Fs}
x\left( z\right) -x\left( \zb\right)  & x^{\prime }\left( z\right)+ x^{\prime }\left( \zb\right) \\
y\left( z\right) -y\left( \zb\right)  & y^{\prime }\left(
z\right)+y^{\prime }\left( \zb\right)
\end{array}%
\right\vert  \\
F_t\left( s,t\right)  &=&-\frac{1}{4}\left\vert
\begin{array}{cc}\label{Ft}
x\left( z\right) -x\left( \zb\right)  & x^{\prime }\left( z\right)-  x^{\prime }\left( \zb\right)\\
y\left( z\right) -y\left( \zb\right)  & y^{\prime }\left(
z\right)-y^{\prime }\left( \zb\right)
\end{array}%
\right\vert
\end{eqnarray}%
We conclude that $F(s,t)$ satisfies $F_s=F_xx_s+F_yy_s$ and
$F_t=F_xx_t+F_yy_t$, where $F_x=-\frac{i}{2}(y(z)-y(\zb))$ and
$F_y=\frac{i}{2}(x(z)-x(\zb))$. Thus
$\nabla(F)=iR(\frac{C(z)-C(\zb)}{2})$. Observe also that, since
$C(z)$ is analytic and real on the real line,
$\nabla(F)=\frac{1}{2}iR(C(z)-C(\zb))$ is real. So $F$ is also
real.

}

\begin{proposition}
 The parameterization
$Q(s,t)=(P(s+it,s-it),F(s,t))$ is isothermal and defines an
improper definite affine sphere $S$. The surface does not depend
on the choice of the analytic parameterization of the curve $C$.
Also $\Omega(s,t)=-i\omega(s,it)$, $A(s,t)=ia(s,it)$ and
$B(s,t)=-ib(s,-it)$. Hence the Pick invariant does not vanish for
strictly convex curves $C$.
\end{proposition}

\proof{ Since $\nabla(F)=iR(\frac{1}{2})(C(z)-C(\zb))$, the same
proof of proposition \ref{GradienteHessiana}, item (3), with
$(z,\zb)$ in place of $(u,v)$ implies that $det(D^2(F))=+1$ and
that $Q(s,t)$ is an isothermal parameterization of a definite
improper affine sphere. If we begin with a different analytic
parameterization $C(\psi(z))$, denote by $F_1$ the corresponding
function. From the above formulas, we have
$\nabla(F_1)(z)=\nabla(F)(\psi(z))\psi'(z)$, which implies that
$F_1(z)=F(\psi(z))$. Thus the surface is independent of the
parameterization.

Since $\Omega(z,\zb)=\frac{i}{4}[C'(z),C'({\zb})]=i[P_z,P_{\zb}]$,
we obtain $\Omega(s,t)=\frac{1}{2}[P_t,P_s]$. On the other hand,
$\omega(s,t)=\frac{1}{2}[p_s,p_t]$. Thus
$\Omega(s,t)=-i\omega(s,it)$.

Finally, $A(s,t)=\frac{i}{4}[C'(s+it),C''(s+it)]=ia(s,it)$ and
$B(s,t)=-\frac{i}{4}[C'(s-it),C''(s-it)]=-ib(s,-it)$.

 }

\subsection{Examples}

Given a curve $C$, assume that we can find a parameterization
$C(s)=(x(s),y(s))$, with $x(z)$ and $y(z)$ analytic functions.

\begin{example}\label{ParabolaOut}

Consider the parabola $C(r)=(r,\frac{r^2}{2})$ of example
\ref{parabola}. Since $X(z,\zb)=\frac{z+\zb}{2}$ and
$Y(s,t)=\frac{z^2+\zb^2}{4}$, we obtain $X(s,t)=s$ and
$Y(s,t)=\frac{s^2-t^2}{2}$. Also, $F(s,t)=if(s,it)=\frac{t^3}{3}$.
So
$$
Q(s,t)=(s,\frac{s^2-t^2}{2},\frac{t^3}{3})
$$
$t>0$, defines a definite improper affine sphere. The area element
of the Blaschke metric is
$\Omega(s,t)=-i\omega(s,it)=\frac{t}{2}$.

\end{example}

\begin{example}\label{CircleOut}
Consider the circle $C(r)=(\cos(r),\sin(r))$ of example
\ref{circle}. Since $X(z,\zb)=\frac{\cos(z)+\cos(\zb)}{2}$ and
$Y(z,\zb)=\frac{\sin(z)+\sin(\zb)}{2}$, we obtain
$X(s,t)=\cos(s)\cosh(t)$ and $Y(s,t)=\sin(s)\cosh(t)$. Also
$F(s,t)=if(s,it)=\frac{1}{2}(\frac{\sinh(2t)}{2}-t)$.
 Hence
$$
Q(s,t)=(\cos(s)\cosh(t),\sin(s)\cosh(t),\frac{1}{2}(\frac{\sinh(2t)}{2}-t))
$$
$t>0$, defines a definite improper affine sphere.  The area
element of the Blaschke metric is
$\Omega(s,t)=-i\omega(s,it)=\frac{\sinh(t)\cosh(t)}{2}$.

\end{example}

\begin{example}\label{HyperbolaOut}

Consider the hyperbola $C(r)=(\exp(r),\exp(-r))$ of example
\ref{hyperbola}. Since $X(z,\zb)=\frac{\exp(z)+\exp(\zb)}{2}$ and
$Y(z,\zb)=\frac{\exp(-z)+\exp(-\zb)}{2}$, we obtain
$X(s,t)=\exp(s)\cos(t)$ and $Y(s,t)=\exp(-s)\cos(t)$. Also,
$F(s,t)=if(s,it)=t-\frac{\sin(2t)}{2}$.
 Hence
$$
Q(s,t)=(\exp(s)\cos(t),\exp(-s)\cos(t),t-\frac{\sin(2t)}{2})
$$
with $t>0$, defines a definite improper affine sphere. The area
element of the Blaschke metric is given by
$\Omega(s,t)=\frac{\sin(2t)}{2}$.
\end{example}

\begin{example}\label{CubicOut}

Consider the cubic $C(r)=(r,r^3)$ of example \ref{cubic}. Since
$X(z,\zb)=\frac{1}{2}(z+\zb)$ and
$Y(z,\zb)=\frac{1}{2}(z^3+\zb^3)$, we obtain $X(s,t)=s$,
$Y(s,t)=s^3-3st^2$. Also $F(s,t)=if(s,it)=2t^3s$.
 Hence
$$
Q(s,t)=(s,s^3-3st^2,2t^3s)
$$
$t>0$, defines a definite improper affine sphere. The area element
of the Blaschke metric is $\Omega(s,t)=-i\omega(s,it)=3st$.
\end{example}

\begin{example}\label{z4}

Consider the curve $C(r)=(r,r^4)$ of example \ref{r4}. Since
$X(z,\zb)=\frac{z+\zb}{2}$ and $Y(z,\zb)=\frac{z^4+\zb^4}{2}$, we
obtain $X(s,t)=s$ and $Y(s,t)=s^4-6s^2t^2+t^4$. Also
$F(s,t)=if(s,it)=4s^2t^3-\frac{4}{5}t^5$. Hence
$$
Q(s,t)=(s,s^4-6s^2t^2+t^4,4s^2t^3-\frac{4}{5}t^5)
$$
$t>0$, defines a definite improper affine sphere. The area element
of the Blaschke metric is $\Omega(s,t)=-i\omega(s,t)=6s^2t-2t^3$.

\end{example}

\subsection{Isothermal tangent lines to the curve}

\subsubsection{Another change of variables}

Consider new coordinates $(u,v)$ defined by $s+t=u$ and $s-t=v$.

\begin{lemma}\label{IsothermalLines}

Fix a point $P(u_0,v_0)$, $u_0>v_0$ and consider the lines
$P(u,v_0)$, $v_0\leq u\leq u_0$ and $P(u_0,v)$, $v_0\leq v\leq
u_0$. These lines touches tangentially the curve at the points
$P(v_0,v_0)$ and $P(u_0,u_0)$, respectively.

\end{lemma}

\proof{ Straightforward calculations shows that
$$
P_u=\frac{1}{4}((1+i)C'(z)+(1-i)C'(\zb))
$$
and
$$
P_v=\frac{1}{4}((1-i)C'(z)+(1+i)C'(\zb)) \ .
$$
Thus, at $P(u_0,u_0)$, $P_v=\frac{C'(u_0)}{2}$. Similarly, at
$P(v_0,v_0)$, $P_u=\frac{C'(v_0)}{2}$. }

\medskip
We shall call the above lines the {\it isothermal tangent lines}.
It is interesting to observe that in the case of a parabola, these
isothermal tangent lines are in fact straight lines (see example
\ref{ParabolaOut}). Also, in examples \ref{ParabolaOut},
\ref{CircleOut}, \ref{HyperbolaOut} and \ref{CubicOut}, the
isothermal lines starting at a point $P(u_0,v_0)$ do not meet $C$
before the tangency points $P(u_0,u_0)$ and $P(v_0,v_0)$. We shall
refer to this property as non-crossing isothermal tangents.
Example \ref{z4} does not have this property, i.e., the isothermal
tangent lines meet $C$ before they arrive at the tangency points.

\subsubsection{Geometric interpretation of $F$ as an outer area distance}

\begin{figure}[htb]
 \centering
 \includegraphics[width=0.35\linewidth]{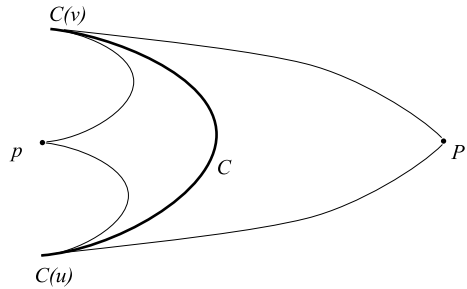}
 \caption{The correspondence between $p$ and $P$.}
\label{Pandp}
\end{figure}

We shall assume from now on that the tangent isothermals are
non-crossing. Considering coordinates $(s,t)$, we can define a
bijection between the inner and the outer parts of $C$, by
corresponding  the points $p=p(s,t)$ and $P=P(s,ti)$. This
correspondence can be defined more geometrically as follows: For
$p\in D$, consider the asymptotic lines that passes through $p$
and denote by $C(u)$ and $C(v)$ the points where they touch
tangentially $C$. Consider then the tangent isothermal lines in
$E$ that touch tangentially $C$ at $C(u)$ and $C(v)$. The
intersection of these lines is $P$ (see figure \ref{Pandp}).

For a point $P=P(u_0,v_0)$, denote by $A(P)$ the region bounded by
the isothermal tangent lines and the part of the curve $C(r)$ with
$v_0\leq r\leq u_0$ (see figure). Next lemma shows that the area
of $A(P)$ is $F(P)/2$. This property justify the name outer area
distance for the function $F$.

\begin{lemma}
Assume that the isothermal tangent lines starting at
$P=P(u_0,v_0)$ are non-crossing. Then the area of $A(P)$ is
$F(P)/2$.

\end{lemma}

\proof{ Straightforward calculations shows that
$[P_u,P_v]=\frac{\Omega}{2}$. Thus the area of $A(P)$ is given by
half of the integral of $\Omega$ over the triangle $T(u_0,v_0)$
whose vertices are $(u_0,v_0)$, $(u_0,u_0)$ and $(v_0,v_0)$. The
corresponding inner region $a(p)$ has area equal to the integral
of $\omega$ over the same triangle. And we know that this area is
equal to $f(p)$. Thus
$$
F(P)=if(p)=i\int_{T(u_0,v_0)}\omega(u,v)dudv=\int_{T(u_0,v_0)}\Omega(u,v)dudv
$$
and so the area of $A(P)$ is $F(P)/2$.

}

\subsubsection{Relation between the indefinite and definite
affine spheres}

\begin{proposition}
Let $b\subset D$ and $B\subset E$ be corresponding regions. Then
$$
\int_{q(b)}j^{1/3}=\int_{Q(B)}J^{1/3}
$$
where the integrals are taken with respect to the Berwald-Blaschke
metric.
\end{proposition}

\proof{
 Remember that $j=ab\omega^{-3}$ and $J=AB\Omega^{-3}$.
 We can assume w.l.o.g. that the curve $C$ is parameterized by
 affine arc-length. Then $j=J=1/16$. Up to these constants,
 the above integrals correspond to the areas, in the $(u,v)$-plane, of the regions
 that represent $b$ and $B$, respectively.
 Since, in $(u,v)$ coordinates, both regions are the same, the proposition is proved.  }

\subsection{Discrete outer area distances to polygons}

Assume that $C$ is a polygon with vertices $P_i=(x_i,y_i)$, $1\leq
i\leq N$. In order to mimic the continuous case, we must consider
that $x$ and $y$ are functions of a variable $u\in Z$ and then
extend these functions to discrete complex analytic functions.

There are several definitions of discrete analytic functions. We
shall adopt a classical one (see \cite{Bobenko-Suris}, ch.5).
Consider $g:Z^2\to R$ and $h:(Z^2)^*\to R$, where $(Z^2)^*$
denotes the dual lattice. $g$ and $h$ are complex conjugates if
$g(u+1,v)-g(u,v)=h(u+1/2,v+1/2)-h(u+1/2,v-1/2)$ and
$g(u,v+1)-g(u,v)=-(h(u+1/2,v+1/2)-h(u-1/2,v+1/2))$, for any
$(u,v)\in Z^2$. These equations are called discrete Cauchy-Riemann
equations, and we say that $g+ih$ is discrete analytic.

Using the above definition, we can extend $x$ and $y$ to discrete
analytic functions $X+ih_1$ and $Y+ih_2$. In fact, these functions
are uniquely defined if we consider that
$h_j(u+1/2,1/2)=-h_j(u+1/2,-1/2)$, $j=1,2$. This condition is
natural for any analytic function that is real on the real line.

As in the continuous case, we must find $F(u,v),(u,v)\in Z^2$ such
that $(h_2,-h_1)=\nabla(F)$. We write
$$
F(u+1,v)-F(u,v)=h_2(u+1/2,v+1/2)\cdot(X(u+1,v)-X(u,v))-h_1(u+1/2,v+1/2)\cdot(Y(u+1,v)-Y(u,v))
$$
and
$$
F(u,v+1)-F(u,v)=h_2(u+1/2,v+1/2)\cdot(X(u,v+1)-X(u,v))-h_1(u+1/2,v+1/2)\cdot(Y(u,v+1)-Y(u,v))
$$
and, using the discrete Cauchy-Riemann equations, we obtain

\begin{eqnarray}
F(u+1,v)-F(u,v)  &=&-\left\vert
\begin{array}{cc}\label{DFu}
h_1(u+1/2,v-1/2)  & h_1(u+1/2,v+1/2) \\
h_2(u+1/2,v-1/2)  & h_2(u+1/2,v+1/2)
\end{array}%
\right\vert  \\
F(u,v+1)-F(u,v)  &=&\left\vert
\begin{array}{cc}\label{DFv}
h_1(u-1/2,v+1/2) & h_1(u+1/2,v+1/2)\\
h_2(u-1/2,v+1/2) & h_2(u+1/2,v+1/2)
\end{array}%
\right\vert
\end{eqnarray}%
In order to simplify notations, we shall denote by $h$ the vector
$(h_1,h_2)$.

\begin{lemma}
There exists $F:Z^2\to R$ satisfying equations (\ref{DFu}) and
(\ref{DFv}).
\end{lemma}

\proof{ One has to show that the discrete derivative $DFu(v)$ of
the right hand of (\ref{DFu}) with respect to $v$ is equal to the
discrete derivative $DFv(u)$ of the right hand of (\ref{DFv}) with
respect to $u$. We have
$$
DFu(v)=[
\begin{array}{cc}
h(u+1/2,v-1/2)+h(u+1/2,v+3/2) & h(u+1/2,v+1/2) \\
\end{array}%
]
$$
and
$$
DFv(u)=[
\begin{array}{cc}
h(u+1/2,v+1/2) & h(u-1/2,v+1/2)+h(u+3/2,v+1/2) \\
\end{array}%
]
$$
Thus $DFv(u)-DFu(v)$ is given by
$$
[
\begin{array}{cc}
h(u+1/2,v+1/2) & h(u-1/2,v+1/2)+h(u+3/2,v+1/2)+h(u+1/2,v-1/2)+h(u+1/2,v+3/2) \\
\end{array}%
]
$$
which is equal to zero, since $h_1$ and $h_2$ are discrete
harmonic. }

\medskip

We shall denote $P=(X,Y)$ and $Q=(X,Y,F)$. Also, define the
co-normal vector by $\nu=(-h_2,h_1,1)=(-\nabla(F),1)$. We can
write
\begin{eqnarray}
Q(u+1,v)-Q(u,v)&=&\nu(u+1/2,v+1/2)\times\nu(u+1/2,v-1/2)\\
Q(u,v+1)-Q(u,v)&=&-\nu(u+1/2,v+1/2)\times\nu(u-1/2,v+1/2).
\end{eqnarray}

It follows directly from these equations that the quadrangles
whose vertices are $Q(u,v),Q(u+1,v),Q(u,v+1)$ and $Q(u+1,v+1)$ are
planar. In fact, each edge is orthogonal to $\nu(u+1/2,v+1/2)$.
Thus $Q(u,v)$ defines a conjugate net (see definition in \cite{
Bobenko-Suris}). Denote by
$\Delta(g)=g(u+1,v)+g(u,v+1)+g(u-1,v)+g(u,v-1)-4g(u,v)$ the
discrete laplacian of $g$. It is  clear that $\Delta(Q)$ is
parallel to the $z$-axis. This follows directly from the fact that
$X$ and $Y$ are discrete harmonic, i.e., $\Delta(X)=\Delta(Y)=0$.
Thus $Q(u,v)$ defines a discrete definite improper affine sphere.

We can also calculate $\Delta(F)$. Denote the vectors
$P(u+1,v)-P(u,v),P(u,v+1)-P(u,v),P(u-1,v)-P(u,v)$ and
$P(u,v-1)-P(u,v)$ by $v_1,v_2,v_3$ and $v_4$, respectively. Since
$P(u,v)$ is discrete harmonic, $v_1+v_2+v_3+v_4=0$. And the area
$A(u,v)$ of the quadrangle whose vertices are
$P(u+1,v),P(u,v+1),P(u-1,v)$ and $P(u,v-1)$ is given by
$[v_3,v_2]+[v_1,v_4]$.

\begin{lemma}
$$
\Delta(F)(u,v)=A(u,v)
$$
\end{lemma}
\proof{ We take as reference the vector $h=h(u+1/2,v+1/2)$. Then
\begin{eqnarray*}
h(u-1/2,v+1/2)&=&h+v_2\\
h(u+1/2,v-1/2)&=&h-v_1\\
h(u-1/2,v-1/2)&=&h+v_2+v_3\\
&=&h-v_1-v_4
\end{eqnarray*}
So
\begin{eqnarray*}
F(u+1,v)-F(u,v)&=&-[h-v_1,h]=[v_1,h]\\
F(u,v+1)-F(u,v)&=&[h+v_2,h]=[v_2,h]\\
F(u-1,v)-F(u,v)&=&[h+v_2+v_3,h+v_2]=[v_3,h]+[v_3,v_2]\\
F(u,v-1)-F(u,v)&=&-[h-v_1-v_4,h-v_1]=[v_4,h]+[v_1,v_4]
\end{eqnarray*}
We conclude that
\begin{eqnarray*}
F(u+1,v)+F(u-1,v)-2F(u,v)=[v_3,v_2]+[v_1,v_4]=A(u,v),
\end{eqnarray*}
thus proving the lemma. }



\begin{figure}[htb]
 \centering
 \includegraphics[width=0.35\linewidth]{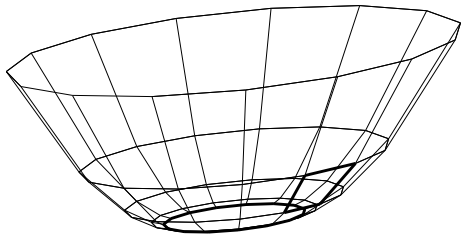}
 \caption{Outer area
distance of a polygon inscribed in a circle: observe the planar
quadrilateral in bold.} \label{DiscreteOutCircle}
\end{figure}

\medskip
One can see in figure \ref{DiscreteOutCircle} the outer area
distance to a polygon inscribed in a circle. In our discrete
construction, it is worthwhile to observe an asymptotic net and a
conjugate net meeting along a curve. Unfortunately we were not
able to give a geometric interpretation of the function $F$ in the
discrete case as we have done in the continuous case.

\begin{figure}[htb]
\centering \fsep \subfigure[ Polygon inscribed in a circle:
$(u,v)$ coordinates inside and $(s,t)$ outside. ] {
\includegraphics[width=.25
\linewidth,clip =false]{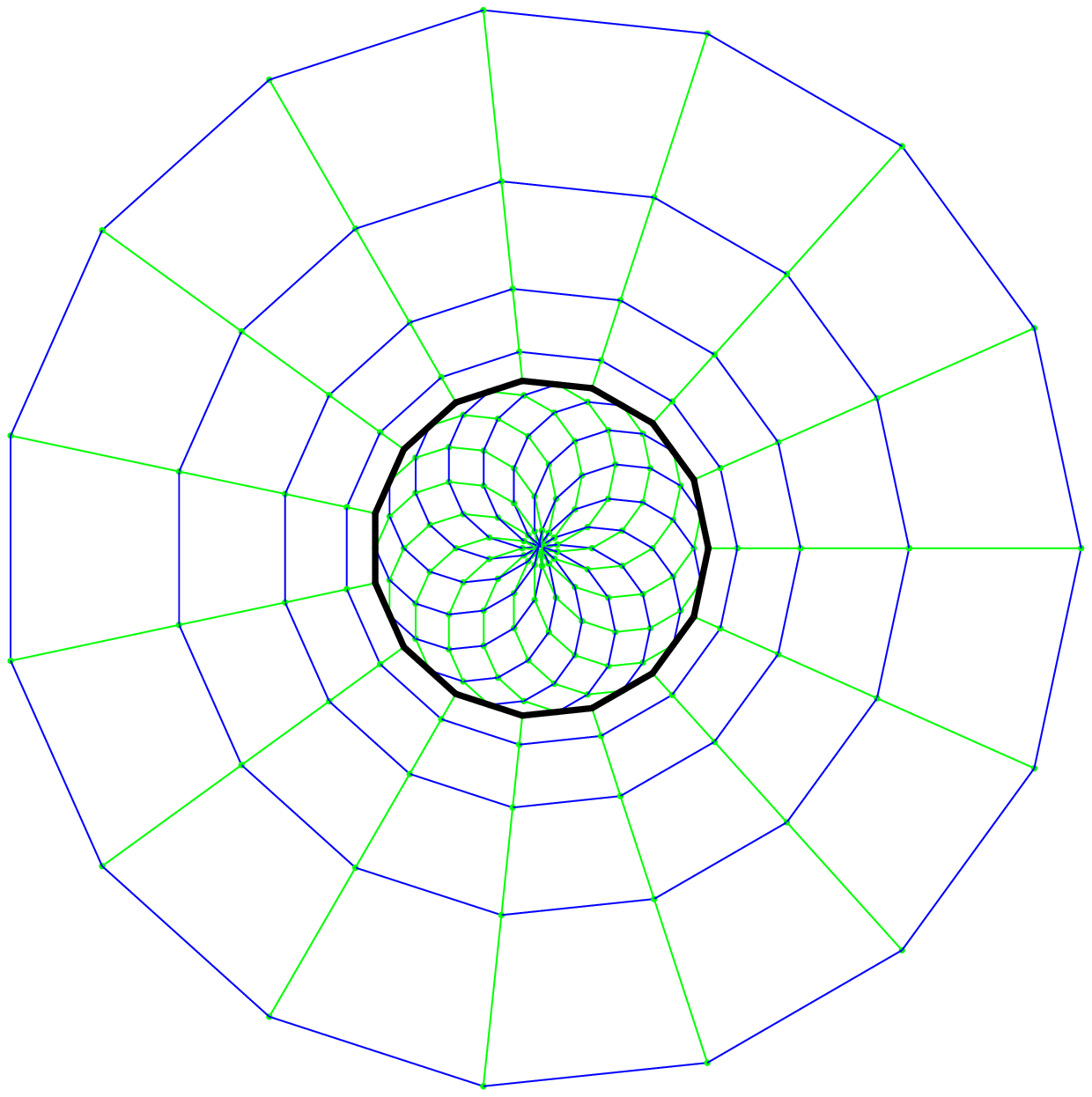}} \fsep\subfigure[ A
closer look. ] {
\includegraphics[width=.25\linewidth,clip
=false]{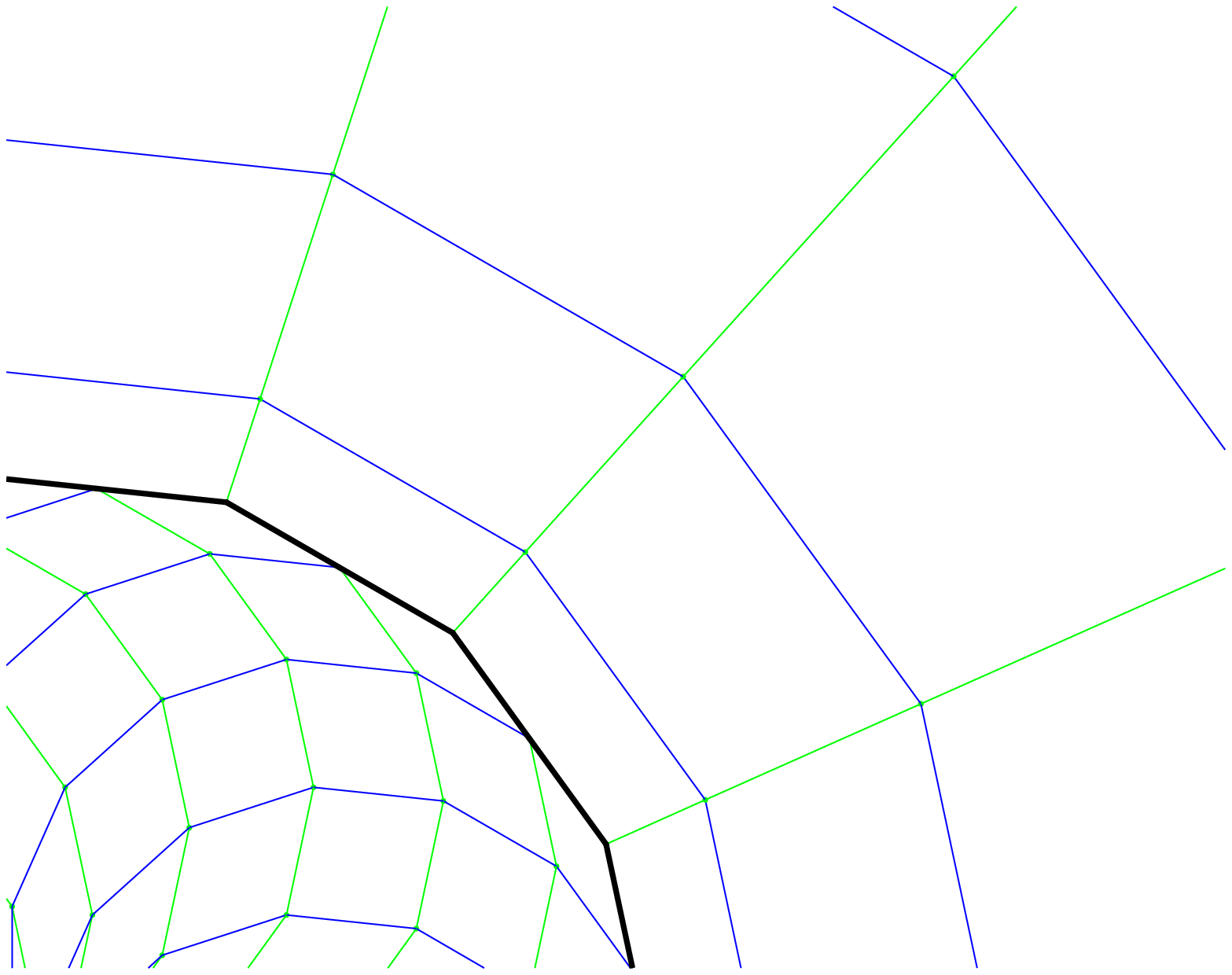}}\fsep\fsep\\
\fsep \subfigure[ Polygon inscribed in a circle:$(u,v)$
coordinates inside and $(u,v)$ outside. ] {
\includegraphics[width=.25
\linewidth,clip =false]{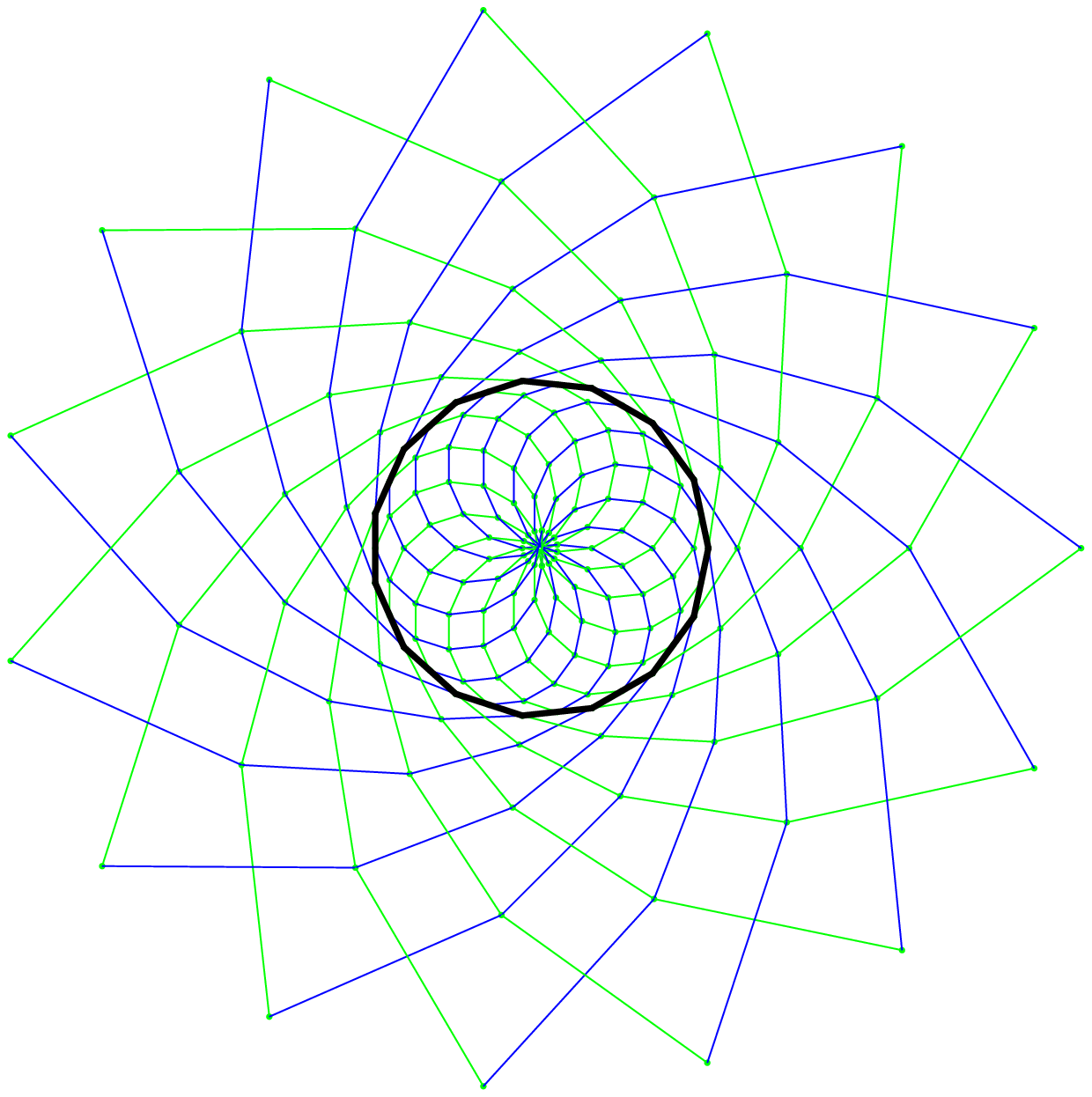}} \fsep\subfigure[ A
closer look. ] {
\includegraphics[width=.25\linewidth,clip
=false]{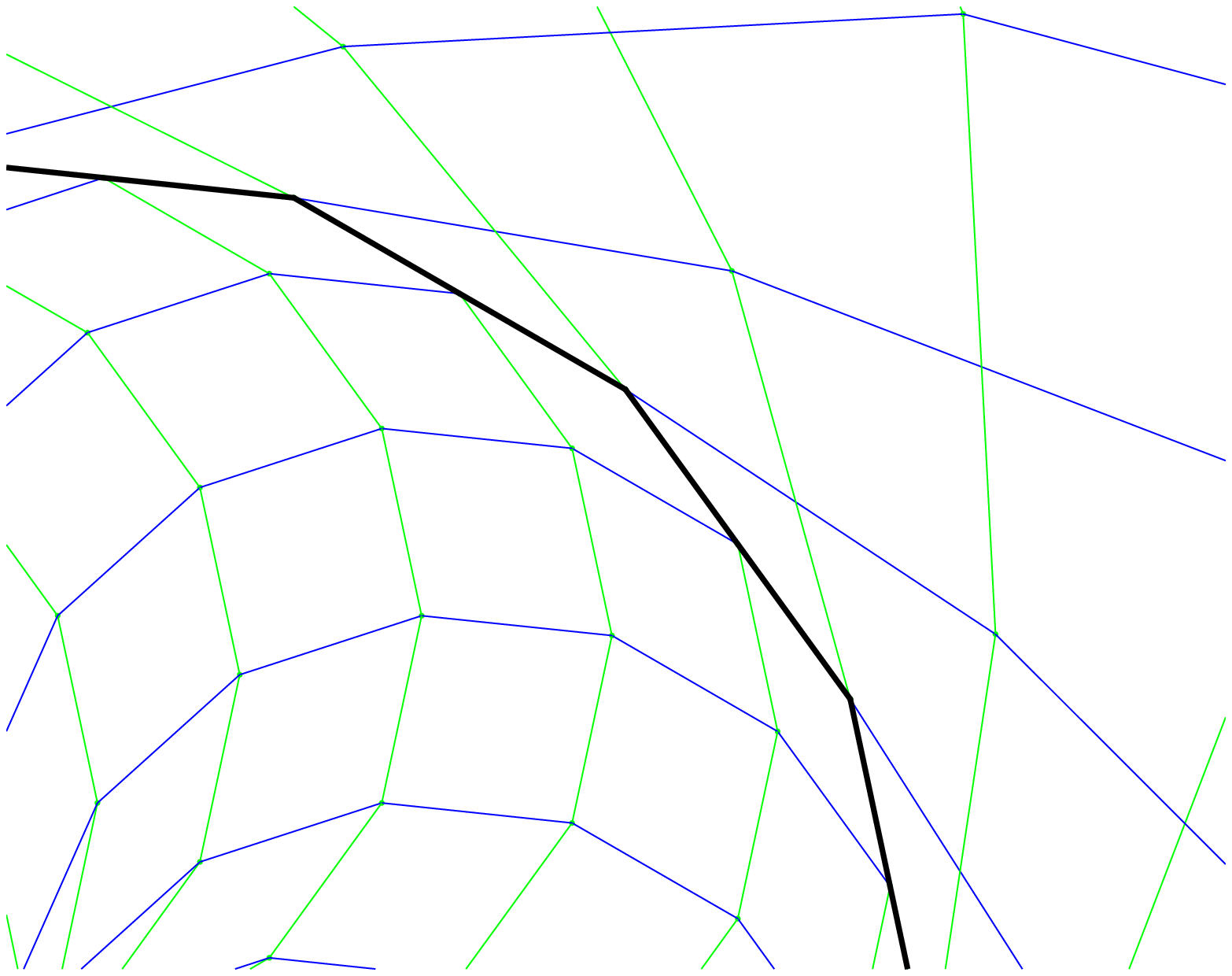}}\fsep\fsep\\
\fsep\subfigure[ Polygon inscribed in a parabola. ]
{\includegraphics[width=.25\linewidth,clip
=false]{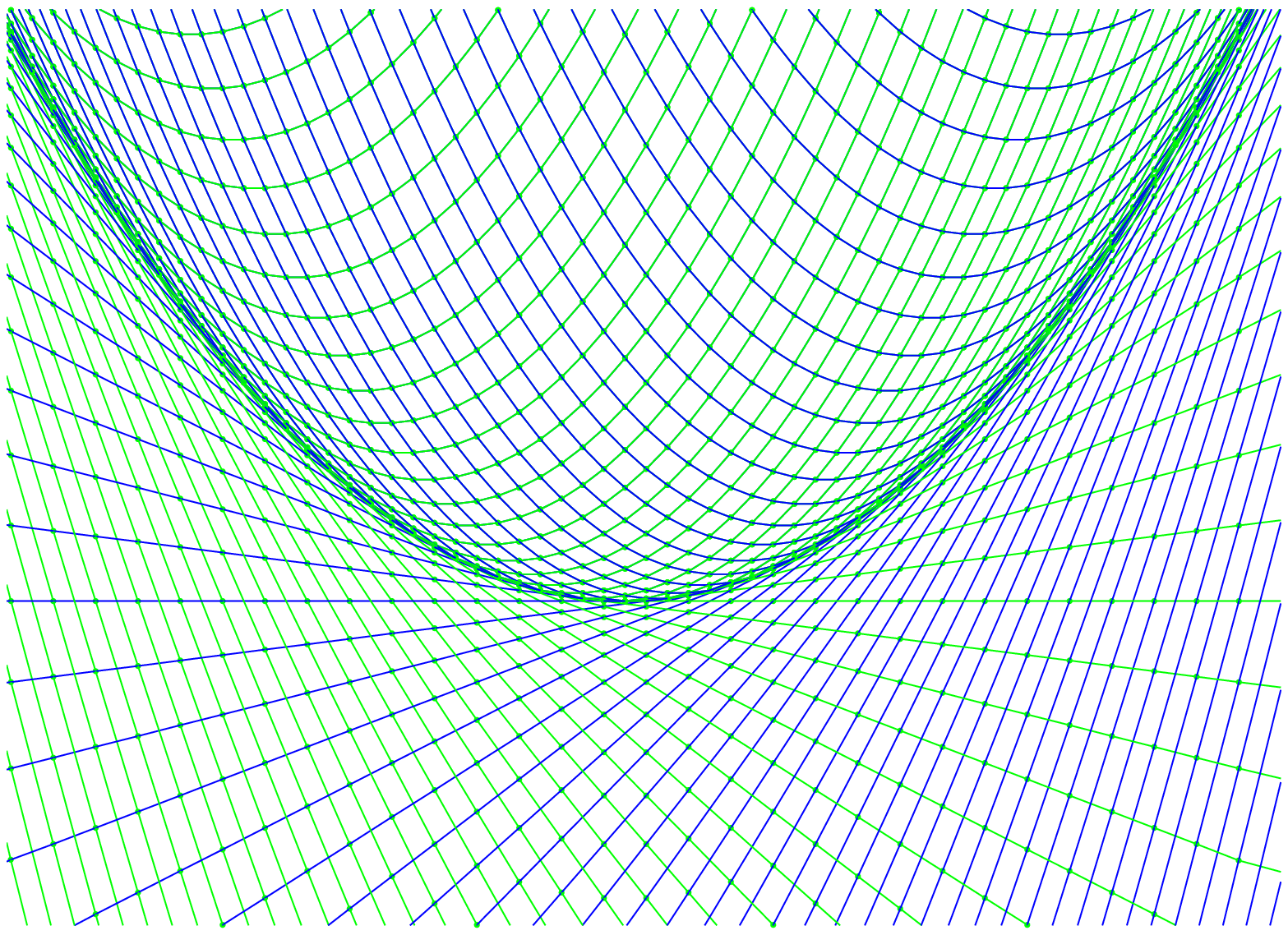}} \fsep\subfigure[ A closer look. ]
{\includegraphics[width=.25\linewidth,clip
=false]{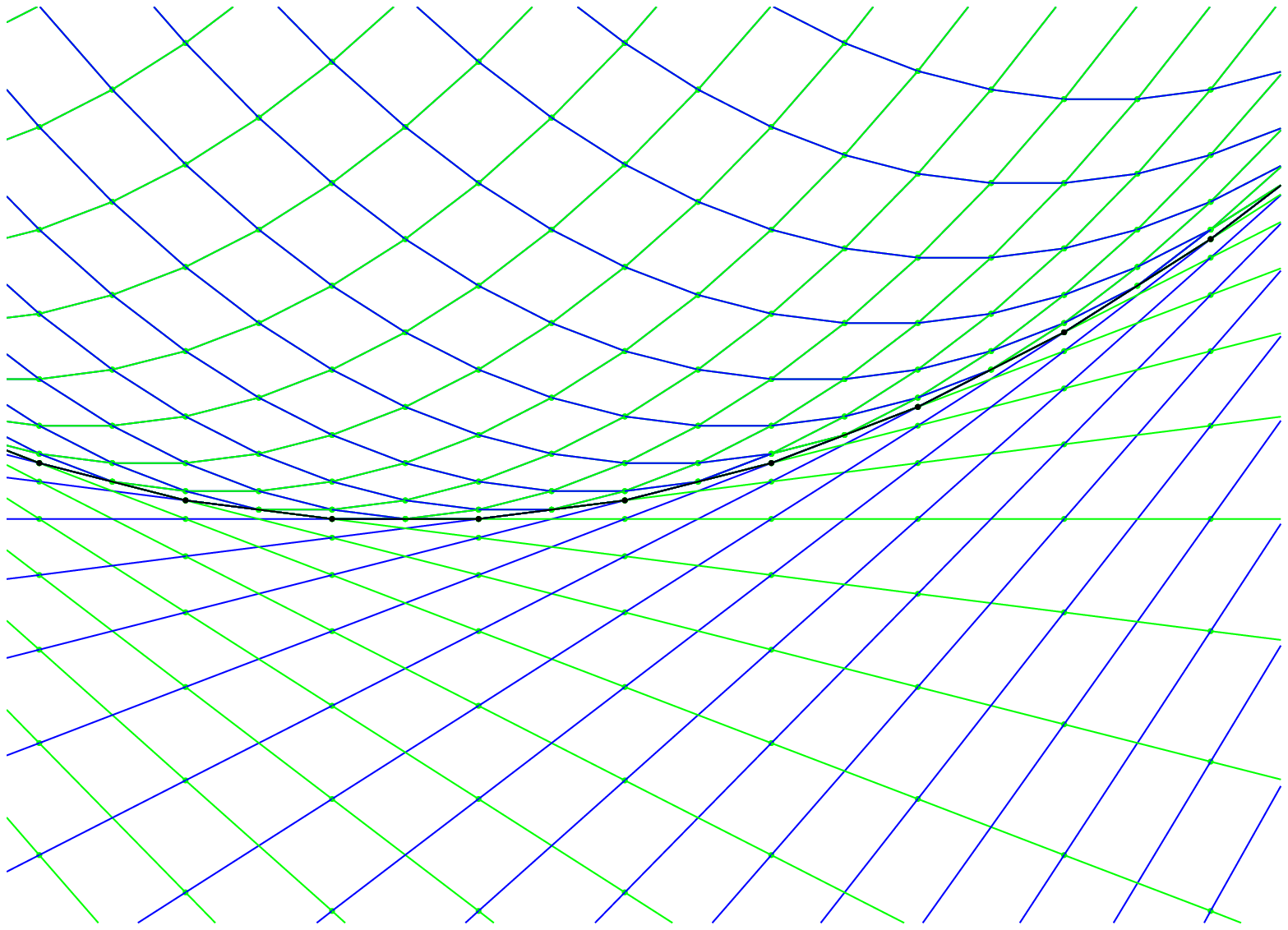}}\fsep \caption{Internal
asymptotic and external isothermal lines of polygons.}
\label{InOutDiscrete2D}
\end{figure}




\section{Conclusions}

We have shown a very close connection between area based distance,
a widely used concept in computer vision, and improper affine
spheres. Both theories may take benefit from this connection.

From the point of view of the theory of area based distances, this
link allow us to propose a new area based distance outside a
convex region and to develop fast algorithms for computing the
inner areas.

From the point of view of the theory of improper affine spheres,
the approach give a very geometrical description of these
surfaces, both in the smooth and in the discrete case.

\smallskip\noindent{\bf Acknowledgements.} The first author thanks CNPq for financial support
through "Projeto Universal", 02/2006.


\bibliography{EsferasAfins}


\end{document}